\documentclass[3p, 11pt]{elsarticle}

\makeatletter
\def\ps@pprintTitle{%
 \let\@oddhead\@empty
 \let\@evenhead\@empty
 \def\@oddfoot{}%
 \let\@evenfoot\@oddfoot}
\makeatother

\usepackage{lipsum}
\usepackage{hyperref}
\usepackage{amssymb}
\usepackage{caption,subcaption}
\usepackage{amsmath}
\usepackage{pgfplots}
\pgfplotsset{compat=1.18}
\usepackage{pgf}
\usepackage{tikz} 
\usepackage{amsfonts}
\usepackage{graphicx}
\usepackage[linesnumbered,ruled,vlined]{algorithm2e}
\usepackage{epstopdf}
\usepackage{MnSymbol}

\usepackage{overpic}
\usepackage{float}
\usepackage{graphicx}
\usepackage[absolute,overlay]{textpos}
\usepackage{booktabs}
\usepackage{relsize}
\usepackage{mathrsfs}

\usepackage[linesnumbered, ruled, vlined]{algorithm2e}
\biboptions{sort&compress}
\newcommand{\bm}[1]{\ensuremath{\mathbf{#1}}}

\newtheorem{definition}{Definition}

\newtheorem{theorem}{\bf Theorem}[section]
\newtheorem{proposition}[theorem]{Proposition}

\begin{document}

\begin{frontmatter}

\title{Implicit Tensor-Train Cross Integration of High-Dimensional Nonlinear PDEs via Fiber-Dependency Elimination}

\author{Behzad Ghahremani and Hessam Babaee$^{*}$}

\address{Department of Mechanical Engineering and Materials Science, University of
Pittsburgh, 3700 O’Hara Street, Pittsburgh, PA, 15213, USA \\ \vspace{2mm}
* Corresponding Author, Email:h.babaee@pitt.edu \vspace{-8mm}}

\begin{abstract} 
Tensor-train (TT) representations have emerged as an effective framework for mitigating the curse of dimensionality in the numerical solution of high-dimensional tensor differential equations. Among existing approaches, TT-cross methods are particularly attractive because they require only pointwise evaluations of the governing equations, naturally accommodate arbitrary nonlinearities, and avoid tangent-space projections and the numerical difficulties associated with nearly singular low-rank factors. However, existing TT-cross rank-truncation methods have been restricted to explicit time integration. Extending TT-cross methods to implicit schemes presents an obstacle: the collocation equations associated with the cross-selected fibers depend on neighboring fibers that are not part of the unknown set. Consequently, the resulting nonlinear system is not closed, preventing the direct application of standard implicit solvers. In this work, we introduce a principled fiber-dependency elimination framework that resolves this obstacle by expressing neighboring fibers as linear combinations of the cross-selected fibers through cross interpolation identities. The resulting formulation produces a closed collocation system while preserving the principal advantages of TT-cross methods. The proposed framework applies to both linear and nonlinear high-dimensional partial differential equations and is naturally combined with Newton iterations and rank adaptivity. Numerical experiments demonstrate rapid convergence of the dependency-elimination iterations, preservation of the temporal accuracy of implicit multistep schemes, and efficient implicit integration of high-dimensional nonlinear problems with full-order discretizations containing up to $10^{55}$ degrees of freedom.
\end{abstract}


\end{frontmatter}

\section{Introduction} \label{sec:Introduction}

Tensor-train (TT) representations \cite{TT11}, originally developed in the tensor-network literature and known in quantum physics as matrix product states \cite{W92}, have emerged as a powerful tool for mitigating the curse of dimensionality in high-dimensional scientific computing. By representing a high-dimensional tensor through a sequence of low-dimensional cores, TT approximations can reduce storage and computational costs from exponential to linear scaling with dimension when low-rank structure is present. These favorable properties have led to successful applications in numerical analysis \cite{K18}, high-dimensional optimization \cite{CLO16}, machine learning \cite{SS16}, and computational fluid dynamics \cite{GLD22}. In particular, TT representations have attracted growing interest for the numerical solution of high-dimensional partial differential equations (PDEs), where conventional discretizations rapidly become intractable due to the exponential growth of the state space.

The use of TT representations for high-dimensional PDEs requires efficient algorithms for evolving low-rank tensors in time. A major advance in this direction was the development of dynamical low-rank (DLR) approximation \cite{KL07,LR13}, which provides a systematic framework for deriving evolution equations directly on low-rank tensor manifolds. In DLR methods, the governing tensor differential equation is projected onto the tangent space of the low-rank manifold, yielding evolution equations for the low-dimensional factors while maintaining the prescribed rank structure. This framework has significantly expanded the applicability of low-rank tensor methods beyond their original quantum-mechanics setting and has led to a variety of algorithms for solving high-dimensional linear and nonlinear PDEs in compressed form.

One important consideration in low-rank time integration is maintaining robustness in the presence of small singular values. Such singular values are unavoidable in practice, as accurate low-rank approximations must resolve progressively finer modes. Within the DLR framework, several strategies have been developed to address this issue. Projector-splitting integrators \cite{CLProjector} avoid explicit inversion of ill-conditioned factors, although they may exhibit instabilities for parabolic problems due to the backward step. Robust basis update and Galerkin (BUG) integrators \cite{CKL22} provide stable first-order schemes, with a second-order midpoint extension proposed in \cite{CEKL24}. Higher-order methods based on rank-truncated time-discrete evolution equations have also been developed for matrix differential equations \cite{ACDynamically,KV19} and tensor differential equations (TDEs) \cite{RDV22}. While these approaches have significantly advanced low-rank time integration, they remain rooted in tangent-space evolution and the associated management of low-rank factors.

Despite substantial progress in low-rank time integration, efficient implicit schemes for nonlinear high-dimensional problems remain comparatively underdeveloped. Within the dynamical low-rank (DLR) framework, several approaches have been proposed to address stiffness, including projected exponential Runge methods \cite{CV23}, enriched-subspace integrators \cite{NQE25,KHSylvester}, BUG formulations \cite{ADRobust,SMPreconditioning}, and Riemannian techniques \cite{SV24}. For TDEs, an implicit integrator on low-rank TT manifolds was introduced in \cite{RV23}. Some of the other broadly used techniques for linear systems in the TT format are alternating least squares (ALS), density matrix renormalization group (DMRG) methods, alternating minimal energy (AMEn), and Krylov methods such as TT-GMRES \cite{HRS12,SD12,DS14,D13}.

Nonlinear problems, however, present an additional challenge. Even when the solution admits an accurate low-rank representation, the nonlinear map of a low-rank tensor generally exhibits substantially higher rank and may become effectively full rank. In high-dimensional settings, this is more than a question of computational cost as the resulting tensor is often too large to be formed or stored explicitly. Consequently, the evaluation of generic nonlinear terms can become a bottleneck, particularly for non-polynomial nonlinearities and within implicit schemes where such evaluations must be performed repeatedly during nonlinear solves. Developing implicit algorithms that can robustly handle nonlinear high-dimensional TDEs without requiring the construction of the full nonlinear tensor therefore remains an important challenge.

An alternative approach is provided by CUR decompositions \cite{GTZ97,OT10}, also referred to as cross or pseudoskeleton approximations. More recently, this framework was employed to efficiently handle nonlinearities in the dynamical low-rank approximation of matrix differential equations \cite{NB23,DPNFB23}. Broadly speaking, cross-based methods can be divided into two classes.  In the first, cross approximation is used within a DLR framework to efficiently evaluate nonlinear terms, either by constructing a cross approximation of the nonlinear tensor or by approximating the projected nonlinear operator \cite{NB23,D25,CKLV25}. While these approaches reduce the cost of nonlinear evaluations, they retain the tangent-space evolution framework. The second class adopts a fundamentally different perspective. Rather than deriving evolution equations for the low-rank factors, it first discretizes the governing equations in time and subsequently constructs a low-rank approximation of the resulting time-discrete solution using cross approximation. Consequently, the low-rank factors are not treated as dynamical variables, and tangent-space projections are avoided altogether \cite{DPNFB23,GBTT24}. This philosophy is closely related to rank-truncated time-discrete evolution methods based on the singular value decomposition (SVD) \cite{RDV22}, as both perform low-rank approximation after time discretization rather than evolving directly on the low-rank manifold. The key distinction is that SVD-based approaches generally require access to the full time-discrete solution prior to truncation. For nonlinear high-dimensional problems, this typically entails forming the nonlinear tensor, which may be prohibitively expensive or impossible to store. In contrast, cross approximation constructs the low-rank solution directly from selected entries, so the nonlinear tensor is never formed explicitly. As a result, arbitrary nonlinear operators can be treated using only pointwise evaluations of the governing equations.

The accuracy of cross approximation methods depends critically on the selection of interpolation fibers. In \cite{NB23,DPNFB23}, the discrete empirical interpolation method (DEIM) \cite{CS10} was adopted for dynamical low-rank approximation, and this strategy has subsequently been employed in several related works; see, for example, \cite{D25,CKLV25}. Nevertheless, the proposed framework is not tied to a particular interpolation strategy. Other sampling techniques, such as QDEIM \cite{DG16} and maximum-volume (MaxVol) methods \cite{GT01}, can also be employed to construct the interpolation fibers.

These applications of cross algorithms have largely been limited to explicit time-integration methods for nonlinear matrix differential equations and tensor differential equations. Extending cross rank-truncation methods to implicit time integration is considerably more challenging because the selected interpolation entries depend on neighboring entries that are themselves unknown, so the resulting nonlinear system is not closed. Consequently, standard implicit solvers cannot be applied directly.

Recent work has shown that dependency-elimination strategies can be constructed in certain matrix settings. For parametric matrix differential equations, the special structure of the governing equations permits one-sided elimination of row dependencies through the evolving low-rank basis \cite{NAB25}. This idea was subsequently generalized to generic matrix differential equations by introducing interpolation coefficient matrices that iteratively eliminate both row and column dependencies \cite{ALB25}. These developments suggest that the low-rank ansatz itself can serve as an algebraic tool for restoring closure of implicit cross discretizations. While implicit cross-based methods have recently been proposed for nonlinear tensor equations in the Tucker format \cite{Appelo25Tucker}, those methods address a fundamentally different tensor representation and employ Anderson acceleration rather than the dependency-elimination strategy developed here.

In this work, we address this challenge by introducing a framework that restores closure of the implicit collocation equations arising in TT-cross rank-truncation methods. To the best of our knowledge, no general framework currently exists for implicit time integration of nonlinear tensor differential equations in TT format based on cross rank truncation; the present work addresses this gap.

More broadly, most existing implicit low-rank solvers are based on Riemannian optimization, tangent-space evolution, or alternating optimization over TT cores \cite{DKS23}. The proposed approach follows a fundamentally different philosophy. Rather than deriving evolution equations for the low-rank factors or solving an optimization problem on the low-rank manifold, it directly solves the implicit collocation equations by restoring their closure through the TT-cross interpolation identities. In this sense, the TT-cross representation is used as an algebraic mechanism for eliminating the unresolved neighboring-fiber dependencies that arise in implicit discretizations.

The main contributions of this paper are summarized as follows:

\begin{itemize}
\item We develop the implicit TT-cross framework for solving high-dimensional linear and nonlinear TDEs. The main contribution is a novel fiber-dependency elimination strategy that transforms the implicit collocation equations into a closed nonlinear system involving only the TT-cross selected fibers.

\item Unlike existing implicit low-rank methods based on Riemannian techniques, tangent-space evolution, or alternating optimization over TT cores, the proposed formulation solves the implicit collocation equations directly in the TT-cross representation. As a result, it requires neither tangent-space projections nor optimization on low-rank manifolds, while preserving the principal advantages of cross rank-truncation methods: the nonlinear tensor is never formed explicitly, only $\mathcal{O}(dnr^2)$ entries of the governing equations are evaluated, and arbitrary nonlinear operators—including non-polynomial nonlinearities such as exponential and fractional functions—are handled naturally.

\item The resulting implicit collocation equations involve at most $\mathcal{O}(nr^2)$ unknowns for each TT core, producing a sequence of low-dimensional nonlinear systems that can be solved efficiently using standard direct or iterative solvers. The proposed formulation therefore retains the computational scalability of explicit TT-cross methods while enabling implicit time integration.


\end{itemize}
The remainder of this paper is organized as follows. Section 2 presents the proposed methodology. Numerical experiments are reported in Section 3, and concluding remarks are given in Section~4.

\section{Methodology}\label{sec:Methodology}

\subsection{Preliminaries}\label{sec:notation}

We briefly summarize the notation used throughout the paper and introduce the tensor train format and the associated low-rank manifolds. Vectors are written in bold lowercase (e.g., $\mathbf x$), matrices in bold uppercase (e.g., $\mathbf X$), and tensors in uppercase (e.g., $X$). Tensor entries are denoted by $X(i_1,\dots,i_d)$. Typewriter font is used for algorithms such as \texttt{SVD} and \texttt{TT-CROSS-DEIM}.  The mode-$k$ unfolding of $X$ is
\[
\mathbf X_{(k)}=
[X(i_k,~i_1\dots i_{k-1}~i_{k+1}\dots i_d)],
\]
which is a matrix of size $n_k\times n_1\dots n_{k-1}n_{k+1}\dots n_d$.

To extract submatrices, we use sets of integer indices. For example, let $\mathcal I=\{i_1^{(\alpha)}\}_{\alpha=1}^r \subset \{1,\dots,n_1\}$. Then $\mathbf X(\mathcal I,:)\in\mathbb R^{r\times n_2}$ contains the selected rows of $\mathbf X\in\mathbb R^{n_1\times n_2}$, $\mathbf X(:,\mathcal J)$ selects columns, and $\mathbf X(\mathcal I,\mathcal J)$ selects the intersection of rows and columns. A colon indicates that all indices of a mode are retained.

Subtensors are handled analogously. For $X\in\mathbb{R}^{n_1\times\dots\times n_d}$ and $\mathcal I=\{i_1^{(\alpha)}\}_{\alpha=1}^r $, the subtensor $X(\mathcal I,:)\in\mathbb R^{r\times n_2\times\dots\times n_d}$ contains the selected indices of the first mode. We also use tuples to denote subsets of multiple modes. For example, if
\[
\mathcal I=\{(i_1^{(\alpha)},i_2^{(\alpha)},i_3^{(\alpha)})\}_{\alpha=1}^r,
\qquad
\mathcal J=\{(i_4^{(\alpha)},\dots,i_d^{(\alpha)})\}_{\alpha=1}^r,
\]
then $X(\mathcal I,:)$ and $X(:,\mathcal J)$ denote subtensors obtained by restricting the corresponding groups of modes. To denote merged modes, we use square brackets. For example,
\[
X([i_1,i_2],i_3,\dots,i_d)
\]
denotes a tensor of size $n_1n_2\times n_3\times\dots\times n_d$, where the first two modes are combined into a single mode.

We denote the mode-$n$ product by $\bigtimes_n$. For
$X\in\mathbb R^{n_1\times\dots\times n_d}$ and
$\mathbf B\in\mathbb R^{m\times n_n}$,
the tensor $X\bigtimes_n\mathbf B$ has size
\[
n_1\times\dots\times n_{n-1}\times m\times n_{n+1}\times\dots\times n_d.
\]
The symbol $\circ$ denotes the Hadamard (elementwise) product.

The Frobenius norm of a tensor is defined as
\begin{equation}\label{eq:FrobNorm}
\|X\|_F=\sqrt{\sum_{i_1=1}^{n_1}\cdots\sum_{i_d=1}^{n_d} X(i_1,\dots,i_d)^2 }.
\end{equation}

For a rank-$r$ truncated singular value decomposition of a matrix
$\mathbf X\in\mathbb R^{m\times n}$, we write
\[
[\mathbf U,\boldsymbol\Sigma,\mathbf Y]
=\texttt{SVD}(\mathbf X,r),
\]
where
$\mathbf U\in\mathbb R^{m\times r}$ contains the left singular vectors,
$\boldsymbol\Sigma\in\mathbb R^{r\times r}$ is the diagonal matrix of singular values, and
$\mathbf Y\in\mathbb R^{n\times r}$ contains the right singular vectors. A tilde symbol ($\sim$) indicates that a factor is discarded. For example,
\[
[\mathbf U,\sim,\sim]
=\texttt{SVD}(\mathbf X,r)
\]
returns only the first $r$ left singular vectors. The order of the tensor modes is rearranged using the \texttt{permute} operation. For example, the dimensions of $X(i_1, i_2, \dots, i_d)$ are reordered as $X(i_d, \dots, i_2, i_1)$ by applying
$X \leftarrow \texttt{permute}(X, [d, d-1, \dots, 1])$.

In the following, we define the tensor train approximation model and the associated low-rank tensor train manifolds.

\begin{definition}[Tensor train format]
A tensor
$\hat X\in\mathbb R^{n_1\times\cdots\times n_d}$
admits a TT representation if its entries can be expressed as
\begin{equation}
\label{eq:TTform}
\hat X(i_1,\dots,i_d)
=
\sum_{\alpha_1=1}^{r_1}
\cdots
\sum_{\alpha_{d-1}=1}^{r_{d-1}}
G_1(1,i_1,\alpha_1)
,G_2(\alpha_1,i_2,\alpha_2)
\cdots
G_d(\alpha_{d-1},i_d,1),
\end{equation}
where
$G_i\in\mathbb R^{r_{i-1}\times n_i\times r_i}$
are the tensor train (TT) cores,
$r_0=r_d=1$, and
$\mathbf r=(r_1,\dots,r_{d-1})$
are the TT ranks.
\end{definition}

\begin{definition}[Low-rank tensor manifolds]
\label{def:Mr}
The low-rank tensor manifold $\mathcal M_r$ is defined as
\[
\mathcal M_r=
{
\hat X\in\mathbb R^{n_1\times\dots\times n_d}
:\operatorname{rank}_{TT}(\hat X)=\mathbf r
},
\]
that is, the set of tensors with fixed TT rank $\mathbf r$. Members of $\mathcal M_r$ are denoted by a hat symbol, e.g., $\hat X$.
\end{definition}

Throughout the paper, we assume $n_i=\mathcal O(n)$ and $r_i=\mathcal O(r)$. Under this assumption, a full tensor contains $\mathcal O(n^d)$ entries, whereas a TT representation requires only $\mathcal O(dnr^2)$ degrees of freedom.

\subsection{TT-CROSS-DEIM for low-rank approximation}
\label{sec:TT-CROSS-DEIM}

We briefly review the \texttt{TT-CROSS-DEIM} algorithm introduced in \cite{GBTT24}, which is closely related to the TT-cross method of \cite{OT10}. The algorithm constructs a TT approximation recursively by applying \texttt{CUR-DEIM} to successive unfoldings of the target tensor. At each step, a small number of fibers are selected through \texttt{DEIM}, and the corresponding TT core is computed from sampled tensor entries. Importantly, the large residual matrices generated during the recursion are never assembled or stored; instead, they are interpreted as subtensors of the original tensor and accessed only through pointwise evaluations. This property is essential in the present setting, where the full tensor may be too large to form explicitly.

Let $X \in \mathbb{R}^{n_1 \times \cdots \times n_d}$ be the target tensor, and let
$\mathbf X_{(1)} \in \mathbb{R}^{n_1 \times n_2 n_3 \cdots n_d}$
denote its mode-1 unfolding. To extract the first TT core, we apply \texttt{CUR-DEIM} to $\mathbf X_{(1)}$. Each column index of $\mathbf X_{(1)}$ corresponds uniquely to a $(d-1)$-tuple of indices. Given a column index set
\[
\mathcal J_1
=
\left\{
\left(i_2^{(\alpha_1)},\dots,i_d^{(\alpha_1)}\right)
\right\}_{\alpha_1=1}^{r_1},
\]
we compute $\mathbf U_1 \in \mathbb{R}^{n_1 \times r_1}$ from
\begin{equation*}
[\mathbf U_1,\sim,\sim]
=
\texttt{SVD}\left(X(:,\mathcal J_1),r_1\right),
\end{equation*}
where $X(:,\mathcal J_1) \in \mathbb{R}^{n_1 \times r_1}$ contains the selected columns of $\mathbf X_{(1)}$. The corresponding \texttt{DEIM} row indices are
\begin{equation*}
\mathcal I_1
=
\texttt{DEIM}(\mathbf U_1),
\qquad
\mathcal I_1
=
\left\{
i_1^{(\alpha_1)}
\right\}_{\alpha_1=1}^{r_1}.
\end{equation*}
The first TT core is obtained from
\begin{equation*}
\mathbf G_1
=
\mathbf U_1
\mathbf U_1(\mathcal I_1,:)^{-1}
\in
\mathbb{R}^{n_1 \times r_1},
\end{equation*}
which is reshaped into a tensor of size $1 \times n_1 \times r_1$ and denoted by $G_1$. The corresponding rank-$r_1$ CUR approximation of $\mathbf X_{(1)}$ is
\[
\hat{\mathbf X}_{(1)}
=
\mathbf G_1 \mathbf R_1,
\qquad
\mathbf R_1
=
\mathbf X_{(1)}(\mathcal I_1,:)
\in
\mathbb{R}^{r_1 \times n_2\cdots n_d}.
\]
Although $\mathbf R_1$ contains $\mathcal{O}(n^{d-1})$ entries, it is never assembled or stored explicitly. The above steps are summarized as
\begin{equation}
[\mathbf G_1,\mathcal I_1]
=
\texttt{CUR-DEIM}(X,\mathcal J_1),
\end{equation}
with \texttt{CUR-DEIM} presented in Algorithm~\ref{alg:CUR-DEIM}. In this setting, $X$ should be interpreted as a function that returns tensor entries at queried indices, rather than as a fully assembled tensor.

The recursion proceeds by applying \texttt{CUR-DEIM} to the residual matrix $\mathbf R_1$. After reshaping, $\mathbf R_1$ is interpreted as a $(d-1)$-dimensional tensor
\[
R_1
\in
\mathbb{R}^{r_1 n_2 \times n_3 \times \cdots \times n_d},
\]
obtained by merging the indices $\alpha_1$ and $i_2$. Using the interpolation indices $\mathcal I_1$, this tensor satisfies
\begin{equation}
\label{eq:aux_1}
R_1(\alpha_1 i_2,i_3,\dots,i_d)
=
X(i_1^{(\alpha_1)},i_2,\dots,i_d),
\qquad
i_1^{(\alpha_1)}\in\mathcal I_1.
\end{equation}
Thus, $R_1$ can be accessed directly through entries of the original tensor $X$ without forming $\mathbf R_1$.

\begin{algorithm}[t]
\fontsize{10pt}{10pt}\selectfont
\SetAlgoLined
\KwIn{$X$: function handle that returns $X(i_1, i_2, \dots, i_m), \quad$  $\mathcal{J}= \{i_2^{(\alpha)}, \dots, i_m^{(\alpha)} \}, \quad \alpha=1, \dots, r$} \vspace{1mm}
\KwOut{$ \mathbf G, ~ \mathcal I$} \vspace{1mm}

$[\mathbf U, \sim, \sim ] = \texttt{SVD}(X(:, \mathcal{J}),r)$ \hspace{1.45cm} $\rhd$ compute the SVD of the selected columns\;

$\mathcal I =\texttt{DEIM}(\mathbf U)$ \hspace{3.50cm} $\rhd$ compute the DEIM indices\;

$\mathbf G =\mathbf U \mathbf U(\mathcal I,:)^{-1}$ \hspace{3.08cm} $\rhd$ compute the column factorized matrix\;

\caption{\texttt{CUR-DEIM} algorithm for unfolded tensors \cite{GBTT24}}
\label{alg:CUR-DEIM}
\end{algorithm}

To extract the second core, we repeat the same steps used for the first core but applied to $R_1$. Let $\mathcal J_2 = \big \{i^{(\alpha_2)}_3, \dots, i^{(\alpha_2)}_d\big \}_{\alpha_2=1}^{r_2}$ be a set of $(d-2)$-tuples, assumed known. Applying \texttt{CUR-DEIM} to $R_1$ gives
\begin{equation}
[\mathbf G_2, \mathcal I_2] = \texttt{CUR-DEIM}(R_1,\mathcal J_2),
\end{equation}
where $\mathbf G_2 \in \mathbb{R}^{r_1 n_2 \times r_2}$. Reshaping $\mathbf G_2$ gives the second TT core
$G_2 \in \mathbb{R}^{r_1 \times n_2 \times r_2}$. The set $\mathcal I_2$ contains $r_2$ DEIM-selected row indices from $\{1,\dots, r_1 n_2\}$. Because of the relation in Eq.\eqref{eq:aux_1}, each row index corresponds uniquely to a pair $\{i_1^{(\alpha_1)}, i_2\}$ where $i_1^{(\alpha_1)} \in \mathcal I_1$ and $i_2 \in \{1,\dots,n_2\}$. Thus,
\begin{equation}
\mathcal I_2 = \big\{ i_1^{(\alpha_2)} , i_2^{(\alpha_2)} \big\}_{\alpha_2=1}^{r_2}, \quad  \mbox{where} \quad {i_1^{(\alpha_2)}} \subset \mathcal I_1.
\end{equation}

These steps are repeated in the same manner to construct all remaining TT cores. The final TT core $G_d$ is obtained by explicitly forming the rows selected by the CUR procedure applied to $R_{d-2}$, a matrix of size $r_{d-2} n_{d-1} \times n_d$. Thus, $\mathbf G_d = R_{d-2}(\mathcal I_{d-1},:) \in \mathbb{R}^{r_{d-1}\times n_d}$, where $\mathcal I_{d-1} \subset \{1,\dots, r_{d-2} n_{d-1}\}$ denotes the DEIM-selected row indices. The matrix $\mathbf G_d$ is then reshaped into the TT core $G_d \in \mathbb{R}^{r_{d-1} \times n_d \times 1}$.

We denote $\mathcal I_1, \mathcal I_2, \dots, \mathcal I_{d-1}$ the left nested indices and $\mathcal J_1, \mathcal J_2, \dots, \mathcal J_{d-1}$ the right nested indices. The procedure described above performs a left-to-right sweep, but a right-to-left sweep is also possible, in which the left indices serve as input and the right indices are produced as output. The earlier steps assumed that suitable sets $\mathcal J_1, \mathcal J_2, \dots, \mathcal J_{d-1}$ are available. In TDEs on low-rank TT manifolds, such near-optimal choices can be taken from the previous time step. In more general settings, however, good right indices are not immediately accessible. A typical example is the construction of a low-rank representation of an initial condition in a TDE. In this setting, the \texttt{\texttt{TT-CROSS-DEIM} (iterative)} procedure is used, where the index sets are refined through iteration. The process is as follows: (i) randomly initialize the right indices $\mathcal J_1, \dots, \mathcal J_{d-1}$; (ii) run a left-to-right sweep to obtain the left indices $\mathcal I_1, \dots, \mathcal I_{d-1}$; (iii) run a right-to-left sweep to update the right indices; (iv) repeat steps (ii)–(iii) until convergence is achieved. The iteration limit can be set based on the convergence of the singular values of the core tensors. Note that $\mathcal{I}$ and $\mathcal{J}$ are nested indices after the first sweep. See \cite{GBTT24} for more details about \texttt{TT-CROSS-DEIM}.

\subsection{Problem setup}
\label{sec:setup}

We consider high-dimensional time-dependent PDEs whose spatial discretization leads to tensor differential equations. Let
\[
u=u(x_1,\dots,x_d,t)
\]
denote a scalar field defined over a $d$-dimensional physical, parametric, or phase-space domain. Throughout this work, we consider TDEs whose spatial operators admit a separable tensor-product representation, in which each operator term acts along a single tensor mode. This class encompasses a broad range of multidimensional PDE discretizations, including all examples considered in this work.

After discretizing each coordinate $x_i$ using $n_i$ degrees of freedom, the semi-discrete solution can be represented as a tensor
\[
X(t)\in\mathbb{R}^{n_1\times n_2\times\cdots\times n_d},
\]
where $X(i_1,\dots,i_d,t)$ denotes the numerical approximation of $u$ at the corresponding grid or basis indices. For clarity of presentation, we first consider the linear TDE
\begin{equation}
\label{eq:linearTDE}
\dot{X}(t)
=
\sum_{i=1}^{d}
X(t)\times_i \mathbf A_i(t)
+
B(t)\circ X(t)
-
C(t),
\end{equation}
supplemented with appropriate initial and boundary conditions. Here,
$\mathbf A_i(t)\in\mathbb{R}^{n_i\times n_i}$ are discrete one-dimensional operators acting along mode $i$,
$B(t),C(t)\in\mathbb{R}^{n_1\times\cdots\times n_d}$ are given tensors,
$\dot X=dX/dt$, and $\circ$ denotes the Hadamard product. In many applications, the matrices $\mathbf A_1,\dots,\mathbf A_d$ are sparse, whereas $B$ and $C$ may be full tensors.

The extension to nonlinear TDEs is discussed in Section~\ref{ExtendNonlinear}; owing to the entrywise nature of cross approximation, the treatment of general nonlinearities does not alter the overall algorithmic structure.

As the resolution increases, the semi-discrete system typically becomes increasingly stiff, making explicit time integration prohibitively expensive. Stiffness may also arise from the governing physics itself, for example through rapidly varying source terms in combustion and chemically reacting flows \cite{JLBC25}.

For clarity of exposition, we first consider the case $B(t)=0$. The extension to $B(t)\neq 0$ is discussed in Section~\ref{ExtendNonlinear} and requires only minor modifications. Applying the implicit Euler scheme to Eq.\eqref{eq:linearTDE} yields
\begin{equation}
\label{implicitTDE2}
X^{k+1}
-
\Delta t
\sum_{i=1}^d
\left[
X^{k+1}\times_i \mathbf A_i
\right]
=
X^k
-
\Delta tC^{k+1},
\end{equation}
where the superscripts $k$ and $k+1$ denote time levels. To simplify the notation, the dependence of the coefficient matrices on time has been suppressed, with the understanding that they are evaluated at $t_{k+1}$.

Equation\eqref{implicitTDE2} may be written compactly as
\begin{equation}
\label{eq:ABform}
\mathcal A(X^{k+1})
=
\mathcal B(X^k),
\end{equation}
where the operators
$\mathcal A,\mathcal B:
\mathbb R^{n_1\times\cdots\times n_d}
\rightarrow
\mathbb R^{n_1\times\cdots\times n_d}$
are defined by
\begin{align*}
\mathcal A(X)
&=
X
-
\Delta t
\sum_{i=1}^d
\left[
X\times_i \mathbf A_i
\right],
\\
\mathcal B(X)
&=
X
-
\Delta t C.
\end{align*}

In high-dimensional settings, directly solving Eq.~\eqref{eq:ABform} is computationally prohibitive. Even storing the solution tensor may be infeasible, since it generally contains $\mathcal O(n^d)$ entries regardless of the sparsity of the coefficient matrices. To overcome this difficulty, we seek low-rank TT approximations of the solution tensors, denoted by $\hat X^k$ and $\hat X^{k+1}$. By representing the solution through a sequence of low-rank cores, the number of unknowns is reduced from $\mathcal O(n^d)$ to $\mathcal O(dnr^2)$. The objective is therefore to solve the time-discrete problem \eqref{eq:ABform} directly on the manifold of low-rank TT tensors, $\mathcal M_r$, while avoiding the construction of full tensors.

\subsection{The closure problem in implicit \texttt{TT-CROSS-DEIM}}\label{sec:closure}
Following the explicit \texttt{TT-CROSS-DEIM} framework of \cite{GBTT24}, we seek the solution at each time step by enforcing the governing equations only on the DEIM-selected fibers. For explicit schemes, this procedure is straightforward because the right-hand side depends only on quantities from the previous time step, which are already known. Consequently, each DEIM-selected fiber can be computed independently.

The situation changes fundamentally for implicit time integration. The unknown solution appears inside the governing operator, and evaluating this operator at a DEIM-selected fiber generally requires neighboring fibers that are not part of the unknown set. As a result, the collocation equations associated with the selected fibers are no longer self-contained.

To illustrate this issue, we consider the implicit Euler discretization. Specifically, we define the residual tensor
\begin{equation}
\label{ResAB}
R
=
\mathcal A(\hat X^{k+1})
-
\mathcal B(\hat X^k),
\end{equation}
and require that $R$ vanish on the DEIM-selected fibers.

Let
\[
\mathcal J_p
=
\left\{
i_{p+1}^{(\alpha_p)},
\dots,
i_d^{(\alpha_p)}
\right\}_{\alpha_p=1}^{r_p},
\qquad
\mathcal I_p
=
\left\{
i_1^{(\alpha_p)},
\dots,
i_p^{(\alpha_p)}
\right\}_{\alpha_p=1}^{r_p},
\]
for $p=1,\dots,d-1$, denote the \texttt{TT-CROSS-DEIM} interpolation indices. The construction of a TT approximation of $\hat X^{k+1}$ requires the computation of the fiber blocks
\[
X^{k+1}(:,\mathcal J_1) \in \mathbb R^{n_1 \times r_1},
\qquad
X^{k+1}([\mathcal I_{p-1},:],\mathcal J_p) \in \mathbb R^{r_1n_2 \times r_2},
\quad
p=2,\dots,d-1,
\]
and
\[
X^{k+1}(\mathcal I_{d-1},:) \in \mathbb R^{r_{d-1} \times n_d}.
\]
We begin by considering the first block, $X^{k+1}(:,\mathcal J_1)$. Enforcing the residual to vanish on the fibers associated with $\mathcal J_1$ yields
\begin{equation}
\label{Core1-1}
X(:,\mathcal J_1)
-
\Delta t
\sum_{i=1}^{d}
[X\times_i \mathbf A_i](:,\mathcal J_1)
=
\hat X^k(:,\mathcal J_1)
-
\Delta t\,C(:,\mathcal J_1),
\end{equation}
where, for notational convenience, we have omitted the superscript \(k+1\), i.e.,
\[
X\equiv X^{k+1},
\qquad
C\equiv C^{k+1}.
\]
The tensor $\hat X^k$ is known from the previous time step and is stored in the TT format.

At first glance, Eq.~\eqref{Core1-1} appears to define a linear system for the unknown fiber block $X(:,\mathcal J_1)$. However, this is generally not the case. When the matrices $\mathbf A_1,\dots,\mathbf A_d$ are diagonal, the equation is closed and can indeed be solved directly for $X(:,\mathcal J_1)$. For general operators, however, the terms
\[
[X\times_i\mathbf A_i](:,\mathcal J_1)
\]
depend on neighboring fibers of $X$ that are not part of the unknown block $X(:,\mathcal J_1)$. Consequently, the selected fibers alone are insufficient to define the collocation system, and Eq.~\eqref{Core1-1} cannot be solved directly.

This observation highlights the important distinction between explicit and implicit cross rank-truncation methods. In explicit schemes, the neighboring fibers are known because they are evaluated from the previous time step. In implicit schemes, however, they are themselves unknown, preventing the construction of a closed collocation system. The remainder of this section develops a dependency-elimination framework that reconstructs these neighboring fibers directly from the DEIM-selected fibers, thereby recovering a closed implicit collocation system.

\subsection{Elimination of neighboring-fiber dependencies}\label{sec:dependency}
\begin{figure}[t]
    \centering   \includegraphics[width=0.6\linewidth]{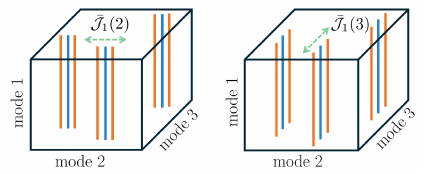} \vspace{-3mm}
    \caption{The union of ${\mathcal{J}}_{1}$ indices and their dependent indices along the second mode (left panel) and third mode (right panel) of $X$.}
    \label{fig:adj-idx}
\end{figure}

 In the following, we present an iterative fiber-dependency elimination procedure that removes these couplings and produces a closed system expressed solely in terms of the unknown block $X(:,\mathcal J_1)$. Let
\[
\bar{\mathcal I}_{p}(i)
=
\left\{
i_1^{(\bar\alpha_p)},
\dots,
i_p^{(\bar\alpha_p)}
\right\}_{\bar\alpha_p=1}^{\bar r_p(i)},
\qquad
p=1,\dots,d-1,
\]
denote the union of $\mathcal I_p$ and the indices coupled to $\mathcal I_p$ through the action of $\mathbf A_i$. Similarly, let
\[
\bar{\mathcal J}_{p}(i)
=
\left\{
i_{p+1}^{(\bar\alpha_p)},
\dots,
i_d^{(\bar\alpha_p)}
\right\}_{\bar\alpha_p=1}^{\bar r_p(i)},
\qquad
p=1,\dots,d-1,
\]
denote the union of $\mathcal J_p$ and the indices coupled to $\mathcal J_p$ through the action of $\mathbf A_i$. Figure~\ref{fig:adj-idx} illustrates the neighboring indices associated with $\mathcal J_1$ along the second and third modes of a three-dimensional tensor.

For $i=2,\dots,d$, the action of the operator $\mathbf A_i$ on the selected fiber block can therefore be written as
\[
[X\times_i \mathbf A_i](:,\mathcal J_1)
=
X(:,\bar{\mathcal J}_1(i))
\times_2
\mathbf A_i
\big(
\mathcal J_1,
\bar{\mathcal J}_1(i)
\big).
\]

The contribution associated with the first mode ($i=1$) requires no additional fibers because all indices of the first mode are already present in the selected block. Consequently,
\[
[X\times_1 \mathbf A_1](:,\mathcal J_1)
=
X(:,\mathcal J_1)\times_1 \mathbf A_1.
\]

Substituting the above expressions into Eq.~\eqref{Core1-1} yields
\begin{multline}
\label{Core1-2}
X(:,\mathcal J_1)
-
\Delta t
\Bigg[
X(:,\mathcal J_1)\times_1 \mathbf A_1
+
\sum_{i=2}^{d}
\Big(
X(:,\bar{\mathcal J}_1(i))
\times_2
\mathbf A_i
(
\mathcal J_1,
\bar{\mathcal J}_1(i)
)
\Big)
\Bigg] 
=
\hat X^k(:,\mathcal J_1)
-
\Delta t\,C(:,\mathcal J_1).
\end{multline}

Eq.~\eqref{Core1-2} makes the source of the difficulty explicit. While the unknown quantity is the DEIM-selected fiber block $X(:,\mathcal J_1)$, the operator evaluations involve neighboring fiber blocks $X(:,\bar{\mathcal J}_1 ( i))$ that are not part of the unknown set. Consequently, the selected fibers alone are insufficient to define a closed collocation system.

The key idea of the proposed methodology is to express these neighboring fibers directly in terms of the DEIM-selected fibers through the \texttt{TT-CROSS-DEIM} representation itself.  As we show below, for the \texttt{TT-CROSS-DEIM} approximation every fiber can be expressed exactly as a linear combination of DEIM-selected fibers. Consequently, the neighboring fibers appearing in Eq.~\eqref{Core1-2} do not constitute additional unknowns; instead, they can be reconstructed directly from the DEIM-selected fibers. This property forms the basis of  the proposed dependency-elimination procedure. In particular, any first-mode fiber of $\hat X$ can be written as
\begin{equation}
\hat{X}(:,i_2,i_3,\dots,i_d)
=
\sum_{\alpha_1=1}^{r_1}
X(
:,
i_2^{(\alpha_1)},
i_3^{(\alpha_1)},
\dots,
i_d^{(\alpha_1)}
)
\,\mathbf z_{r_1}(\alpha_1),
\label{eq:deim_dep1}
\end{equation}
where $\mathbf z_{r_1}$ denotes a vector of interpolation coefficients.

Similarly, fibers associated with an intermediate mode satisfy
\begin{multline}
\hat{X}(i_1,\dots,i_{p-1},:,i_{p+1},\dots,i_d)
 =
\sum_{\alpha_{p-1}=1}^{r_{p-1}}
\sum_{\alpha_p=1}^{r_p}
\mathbf z_{l_p}(\alpha_{p-1})
X(
i_1^{(\alpha_{p-1})},
\dots,
i_{p-1}^{(\alpha_{p-1})},
:,
i_{p+1}^{(\alpha_p)},
\dots,
i_d^{(\alpha_p)}
)
\mathbf z_{r_p}(\alpha_p),
\label{eq:deim_dep2}
\end{multline}
where $\mathbf z_{l_p}$ and $\mathbf z_{r_p}$ are interpolation coefficients associated with the left and right TT contractions, respectively.

Finally, fibers along the last mode satisfy
\begin{equation}
\hat{X}(i_1,\dots,i_{d-1},:)
=
\sum_{\alpha_{d-1}=1}^{r_{d-1}}
\mathbf z_{l_d}(\alpha_{d-1})
X(
i_1^{(\alpha_{d-1})},
\dots,
i_{d-1}^{(\alpha_{d-1})},
:
),
\label{eq:deim_dep3}
\end{equation}
where $\mathbf z_{l_d}$ denotes the corresponding interpolation coefficients.

The derivation of Eqs.~\eqref{eq:deim_dep1}--\eqref{eq:deim_dep3} follows directly from the \texttt{TT-CROSS-DEIM} construction and is provided in \ref{appendix:deimrelations}. These relations show that the neighboring fibers appearing in Eq.~\eqref{Core1-2} can, in principle, be represented through the DEIM-selected fibers. The remaining task is to compute the interpolation coefficients $\mathbf z_{l_p}$ and $\mathbf z_{r_p}$ in a way that is consistent with the \texttt{TT-CROSS-DEIM} representation.  Next, we present a coefficient-update procedure and show how it converts Eq.~\eqref{Core1-2} into a closed nonlinear collocation system.
\begin{figure}[t]
    \centering
        \includegraphics[width=.9\linewidth]{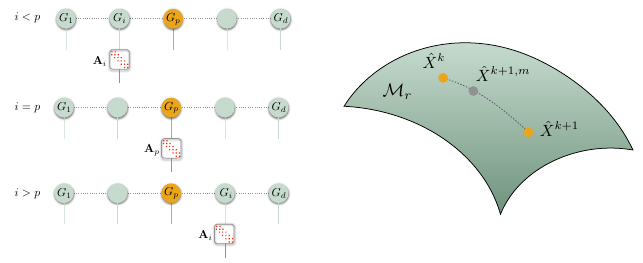}
    \caption{Illustration of the proposed implicit \texttt{TT-CROSS-DEIM} algorithm. Left: Only operator terms associated with modes other than the active TT core produce neighboring-fiber dependencies. Right: Each dependency-elimination iteration reconstructs a valid fixed-rank \texttt{TT-CROSS-DEIM} approximation, so all intermediate iterates remain on the TT manifold $\mathcal M_r$ throughout the nonlinear solve.}
    \label{fig:overview}
\end{figure}

\subsection{Closed implicit collocation equations}
\label{implictTTDEIM}

The interpolation relations developed in the previous section provide a mechanism for eliminating the neighboring-fiber dependencies that prevent the implicit collocation equations from closing.Throughout this section, $\mathbf Z_{l_p}$ and $\mathbf Z_{r_p}$ denote the interpolation coefficient matrices associated with the left and right neighboring fibers, respectively, while $\mathbf Z$ is used as a generic notation when the distinction is not important.  Once the coefficient matrices $\mathbf Z$ have been computed, all dependent fibers can be expressed directly in terms of the DEIM-selected fibers. Consequently, the implicit residual equations reduce to a sequence of closed systems for the \texttt{TT-CROSS-DEIM} subtensors.

We now derive the resulting equations for the first, middle, and last TT cores. 

\subsubsection{First TT core}

For the first TT core, the closure relation
\begin{equation}
X(:,\bar{\mathcal J}_1(i))
=
X(:,\mathcal J_1)\,
\mathbf Z_{r_1}(i),
\qquad
i=2,\dots,d,
\label{eq:closure_first}
\end{equation}
where $\mathbf Z_{r_1}(i) \in \mathbb R^{r_1 \times \bar{r}_1(i)}$ eliminates all neighboring-fiber dependencies. Substituting Eq.~\eqref{eq:closure_first} into the residual-collocation equations yields
\begin{multline}
\label{Core1-3}
    X(:, \mathcal{J}_1)
    - \Delta t
    \Bigg[
    X(:, \mathcal{J}_1)\times_1 \mathbf A_1
    +
    \sum_{i=2}^{d}
    X(:, \mathcal{J}_1)
    \times_2
    \Big(\mathbf A_i
    (\mathcal{J}_1,\bar{\mathcal J}_{1}(i)) \left(\mathbf Z_{r_1}(i) \right)^{\top} \Big)
    \Bigg]
    =
    \hat{X}^k(:, \mathcal{J}_1)
    -
    \Delta t\, C(:, \mathcal{J}_1).
\end{multline}

For fixed coefficient matrices $\mathbf Z_{r_1}(i)$, Eq.~\eqref{Core1-3} is a linear matrix equation for the unknown block $X(:,\mathcal J_1)\in\mathbb R^{n_1\times r_1}$.
Depending on the operator structure, the resulting problem is either a Sylvester equation or a generalized linear matrix equation and can be solved efficiently using direct or iterative methods. In the specific form of TDE chosen for this derivation, Eq. (\ref{Core1-3}) results in a Sylvester equation of the form  of 
\begin{equation}\label{eq:sylv1}
\bm L_1 \bm X + \bm X \bm R_1^{\top} = \bm F_1, 
\end{equation}
where:
\begin{align*}
    & \bm X = X (:, \mathcal{J}_1),\\
    & \bm L_1 = \bm I_{n_1} - \Delta t \: \bm A_1,\\
    & \bm R_1 = - \Delta t \: \sum_{i=2}^{d} \mathbf A_i
    (\mathcal{J}_1,\bar{\mathcal J}_{1}(i)) \left(\mathbf Z_{r_1}(i) \right)^{\top}  ,\\
    & \bm F_1 = \hat X^{k}(:, \mathcal{J}_1) - \Delta t\, C(:, \mathcal J_1).    
\end{align*}
where $ 
\mathbf L_1\in\mathbb R^{n_1\times n_1}$ and $\mathbf R_1\in\mathbb R^{r_1\times r_1}$. The above equation is a low-dimensional Sylvester equation, or more generally a linear matrix equation, which can be solved efficiently using either direct methods, such as the Bartels–Stewart algorithm \cite{BS72}, or iterative methods, including alternating direction implicit (ADI) iterations \cite{BLP08} and Krylov subspace methods \cite{S90}, when the system size is sufficiently large. We refer the reader to \cite{Si16} for a comprehensive review of numerical methods for large-scale linear matrix equations.

\subsubsection{Middle TT cores}

For a middle TT core, both left and right neighboring-fiber dependencies may arise. Applying the closure relations
\begin{equation}
X(\bar{\mathcal I}_{p-1}(i),:,\mathcal J_p)
=
X(\mathcal I_{p-1},:,\mathcal J_p)
\times_1
\mathbf Z_{l_p}(i),
\qquad
i=1,\dots,p-1,
\label{eq:closure_left}
\end{equation}
and
\begin{equation}
X(\mathcal I_{p-1},:,\bar{\mathcal J}_p(i))
=
X(\mathcal I_{p-1},:,\mathcal J_p)
\times_3
\left(\mathbf Z_{r_p}(i)\right)^{\top},
\qquad
i=p+1,\dots,d,
\label{eq:closure_right}
\end{equation}
where $\mathbf Z_{l_p}(i)
\in
\mathbb R^{\bar r_{p-1}(i)\times r_{p-1}}$ and $\mathbf Z_{r_p}(i)
\in
\mathbb R^{r_p\times \bar r_p(i)}$, yields the closed equation
\begin{multline}
\label{Core2-3}
X(\mathcal I_{p-1},:,\mathcal J_p)
-
\Delta t
\Bigg[
\sum_{i=1}^{p-1}
X(\mathcal I_{p-1},:,\mathcal J_p)
\times_1
\Big(
\mathbf A_i(\mathcal I_{p-1},\bar{\mathcal I}_{p-1}(i))
\mathbf Z_{l_p}(i)
\Big)
\\
+
X(\mathcal I_{p-1},:,\mathcal J_p)
\times_2
\mathbf A_p
+
\sum_{i=p+1}^{d}
X(\mathcal I_{p-1},:,\mathcal J_p)
\times_3
\Big(
\mathbf A_i(\mathcal J_p,\bar{\mathcal J}_p(i))
\left(\mathbf Z_{r_p}(i)\right)^{\top}
\Big)
\Bigg]
\\
=
\hat X^k(\mathcal I_{p-1},:,\mathcal J_p)
-
\Delta t\,
C(\mathcal I_{p-1},:,\mathcal J_p).
\end{multline}
The unknown subtensor $X(\mathcal I_{p-1},:,\mathcal J_p)
\in
\mathbb R^{r_{p-1}\times n_p\times r_p}$
contains only $\mathcal O(r_{p-1}n_pr_p)$ entries. After reshaping the first two modes, Eq.~\eqref{Core2-3} becomes a low-dimensional Sylvester equation of form 
\begin{equation}
\label{eq:sylv2}
\bm L_p \bm X + \bm X \bm R_p^{\top} =  \bm F_p,
\end{equation}
where
\begin{align*}
    & \bm X =  X \bigl ( [\mathcal {I}_{p-1}, :] , \mathcal{J}_p \bigr ), \\
    & \bm L_p = \bm I_{n_p r_{p-1}} - {\Delta t}     \sum_{i=1}^{p-1} \Big(\mathbf A_i(\mathcal I_{p-1},\bar{\mathcal I}_{p-1}(i)) \mathbf Z_{l_p}(i) \Big) \otimes \mathbf{I}_{n_p}  \: + \:  \bm I_{r_{p-1}} \otimes \bm A_{p},    \\ 
    & \bm R_p = - \Delta t \:  \sum_{i=p+1}^{d} \mathbf A_i
    (\mathcal{J}_p,\bar{\mathcal J}_{p}(i)) \left(\mathbf Z_{r_p}(i) \right)^{\top}, \\
    & \bm F_p = \hat X^k \big ([\mathcal I_{p-1},:],\mathcal J_p \big) - \Delta t\, C \big ([\mathcal I_{p-1},:],\mathcal J_p \big). 
\end{align*}
where  $\mathbf L_p\in\mathbb R^{r_{p-1}n_p\times r_{p-1}n_p}$, and $
\mathbf R_p\in\mathbb R^{r_p\times r_p}$. Figure \ref{fig:overview} (left side) shows the schematics of the terms for $i<p$, $i=p$ and $i>p$. 

\subsubsection{Last TT core}

For the final TT core, the closure relation
\begin{equation}
X(\bar{\mathcal I}_{d-1}(i),:)
=
\mathbf Z_{l_d}(i)
X(\mathcal I_{d-1},:),
\qquad
i=1,\dots,d-1,
\label{eq:closure_last}
\end{equation}
where $\mathbf Z_{l_d} \in \mathbb R^{\bar{r}_{d-1}\times r_{d-1}}$,
eliminates the remaining neighboring-fiber dependencies. Substituting Eq.~\eqref{eq:closure_last} into the collocation equations yields
\begin{multline}
\label{Core3-3}
X(\mathcal I_{d-1},:)
-
\Delta t
\Bigg[
\sum_{i=1}^{d-1}
X(\mathcal I_{d-1},:) \times_1 \Big(
\mathbf A_i(\mathcal I_{d-1},\bar{\mathcal I}_{d-1}(i))
\mathbf Z_{l_d}(i)
\Big)
\\
+
X(\mathcal I_{d-1},:)
\times_2
\mathbf A_d
\Bigg]
=
\hat X^k(\mathcal I_{d-1},:)
-
\Delta t\,
C(\mathcal I_{d-1},:).
\end{multline}

The resulting problem is again a low-dimensional linear matrix equation for $ X(\mathcal I_{d-1},:)
\in
\mathbb R^{r_{d-1}\times n_d},$
which can be solved efficiently using the same techniques employed for the previous cores. Thus, the above equation may also result in a Sylvester equation of form 
\begin{equation}\label{eq:sylv3}
    \bm L_d \bm X + \bm X \bm R_d^{\top} =  \bm F_d
\end{equation} where:
\begin{align*}
    & \bm X = X (\mathcal{I}_{d-1}, :) \\
    & \bm L_d = \bm I_{r_{d-1}} - {\Delta t} \:    \sum_{i=1}^{d-1} \mathbf A_i(\mathcal I_{d-1},\bar{\mathcal I}_{d-1}(i)) \mathbf Z_{l_d}(i)   \\   
    & \bm R_d = - \: \Delta t \: \bm A_{d}  \\
    & \bm F_d =  \hat X^k(\mathcal I_{d-1},:) - \Delta t\, C(\mathcal I_{d-1},:). 
\end{align*}
where $\mathbf L_d\in\mathbb R^{r_{d-1}\times r_{d-1}}$ and $
\mathbf R_d\in\mathbb R^{n_d\times n_d}$.

The above derivation exploits the separable operator structure of Eq.~\eqref{eq:linearTDE}, in which each operator term acts on a single tensor mode. Consequently, every operator contribution introduces neighboring-fiber dependencies on only one side of an intermediate TT core, leading to the one-sided closure relations derived above. The dependency-elimination framework itself is not restricted to this setting. For more general nonseparable operators, such as mixed-derivative terms, individual operator terms may simultaneously couple both the left and right interpolation sets of a TT core, resulting in two-sided closure relations involving both $\mathbf Z_{l_p}$ and $\mathbf Z_{r_p}$. In that case, the interpolation coefficient matrices must satisfy a coupled nonlinear system and may be computed, for example, through an alternating fixed-point iteration in which the left and right coefficient matrices are updated sequentially while holding the other fixed. The development and analysis of such extensions are left for future work.

The above equations assume that the interpolation coefficient matrices
$\mathbf Z$ are available.
Next, we develop a fixed-point iteration for updating
these matrices in a manner consistent with the \texttt{TT-CROSS-DEIM} representation.

\subsection{Computation of interpolation coefficients}\label{sec:interp_coef}
\label{sec:Zupdate}
The closed collocation systems derived in the previous subsection depend on the interpolation coefficient matrices
$\mathbf Z_{r_1}(i)$,
$\mathbf Z_{l_p}(i)$,
$\mathbf Z_{r_p}(i)$,
and
$\mathbf Z_{l_d}(i)$.
Since these coefficients depend on the current TT approximation, they cannot be determined a priori. We utilize a fixed-point iteration in which the TT cores and the interpolation coefficient matrices are updated alternately until convergence.

The derivation relies on a fundamental interpolation property of the \texttt{TT-CROSS-DEIM} construction. To state this property, we first introduce the contracted TT-core vectors
\begin{align*}
\mathbf g_{1\rightarrow p}^{(\alpha_p)}
&=
G_1(1,i_1^{(\alpha_p)},:)
G_2(:,i_2^{(\alpha_p)},:)
\cdots
G_p(:,i_p^{(\alpha_p)},:),
\\
\mathbf g_{p+1\rightarrow d}^{(\alpha_p)}
&=
G_{p+1}(:,i_{p+1}^{(\alpha_p)},:)
\cdots
G_d(:,i_d^{(\alpha_p)},1),
\end{align*}
where
$\mathbf g_{1\rightarrow p}^{(\alpha_p)}\in\mathbb R^{1\times r_p}$
are row vectors and
$\mathbf g_{p+1\rightarrow d}^{(\alpha_p)}\in\mathbb R^{r_p\times1}$
are column vectors. Similarly, for the neighboring interpolation indices we define
\begin{align*}
\bar{\mathbf g}_{1\rightarrow p}^{(\bar\alpha_p)}
&=
G_1(1,i_1^{(\bar\alpha_p)},:)
\cdots
G_p(:,i_p^{(\bar\alpha_p)},:),
\\
\bar{\mathbf g}_{p+1\rightarrow d}^{(\bar\alpha_p)}
&=
G_{p+1}(:,i_{p+1}^{(\bar\alpha_p)},:)
\cdots
G_d(:,i_d^{(\bar\alpha_p)},1),
\end{align*}
where
$\bar{\mathbf g}_{1\rightarrow p}^{(\bar\alpha_p)}\in\mathbb R^{1\times r_p}$
and
$\bar{\mathbf g}_{p+1\rightarrow d}^{(\bar\alpha_p)}\in\mathbb R^{r_p\times1}$. Stacking these vectors gives
\[
\mathbf G_{1\rightarrow p}
=
\begin{bmatrix}
\mathbf g_{1\rightarrow p}^{(1)}
\\
\vdots
\\
\mathbf g_{1\rightarrow p}^{(r_p)}
\end{bmatrix},
\qquad
\mathbf G_{p+1\rightarrow d}
=
\begin{bmatrix}
\mathbf g_{p+1\rightarrow d}^{(1)}
&
\cdots
&
\mathbf g_{p+1\rightarrow d}^{(r_p)}
\end{bmatrix},
\]
where $\mathbf G_{1\rightarrow p} \in \mathbb R^{r_p \times r_p}$ and $\mathbf G_{p+1\rightarrow d} \in \mathbb R^{r_p \times r_p}$. 
Similarly, $\bar{\mathbf G}_{1\rightarrow p} \in \mathbb R^{\bar{r}_p \times r_p}$ and $\bar{\mathbf G}_{p+1\rightarrow d} \in \mathbb R^{r_p \times \bar{r}_p}$ are obtained by stacking  $\bar{\mathbf g}_{1\rightarrow p}^{(\bar\alpha_p)}$ and $\bar{\mathbf g}_{p+1\rightarrow d}^{(\bar\alpha_p)}$ vectors for $\bar\alpha_p = 1,\dots, \bar{r}_p$.  
The following proposition is the key algebraic property underlying the dependency-elimination procedure.
\begin{proposition}[Interpolation identities]
\label{prop:Gidentity}
The contracted TT-core matrices evaluated at the DEIM interpolation indices satisfy
\begin{equation}
\mathbf G_{1\rightarrow p}
=
\mathbf I_{r_p},
\label{eq:Gidentity_left}
\end{equation}
for every forward \texttt{TT-CROSS-DEIM} sweep. Likewise, after the reverse \texttt{TT-CROSS-DEIM} sweep,
\begin{equation}
\mathbf G_{p+1\rightarrow d}
=
\mathbf I_{r_p}.
\label{eq:Gidentity_right}
\end{equation}
\end{proposition}

The proof is given in \ref{appendix:identity}. As shown below, the above proposition allows the interpolation coefficient matrices to be obtained directly from neighboring TT contractions without solving additional linear systems.

We derive the interpolation coefficient matrices by substituting the TT representation into the closure relations
Eqs.~\eqref{eq:closure_first}, \eqref{eq:closure_left}, \eqref{eq:closure_right}, and \eqref{eq:closure_last}. Substituting the TT representation into both sides of Eq.~\eqref{eq:closure_first} gives, columnwise,
\[
G_1\,\bar{\mathbf g}_{2\rightarrow d}^{(\bar\alpha)}
=
G_1
\sum_{\alpha=1}^{r_1}
\mathbf g_{2\rightarrow d}^{(\alpha)}
\mathbf Z_{r_1}(i)(\alpha,\bar\alpha),
\qquad
\bar\alpha=1,\dots,\bar r_1(i),
\]
where products involving \(G_1\) are understood after reshaping the first TT core into the matrix $\mathbf G_1\in\mathbb R^{n_1\times r_1}$. Stacking these relations over all \(\bar\alpha\) yields the compact form
\[
\mathbf G_1\,
\bar{\mathbf G}_{2\rightarrow d}(i)
=
\mathbf G_1\,
\mathbf G_{2\rightarrow d}\,
\mathbf Z_{r_1}(i).
\]

By construction, $\mathbf G_1(\mathcal I_1,:)=\mathbf I_{r_1},$
implying that \(\mathbf G_1\) has full column rank and therefore defines an injective linear map. Consequently, the common left factor \(\mathbf G_1\) can be eliminated as shown below
\[
\bar{\mathbf G}_{2\rightarrow d}(i)
=
\mathbf G_{2\rightarrow d}
\mathbf Z_{r_1}(i).
\]

According to Proposition \ref{prop:Gidentity}, after the reverse \texttt{TT-CROSS-DEIM} sweep, the interpolation identity $\mathbf G_{2\rightarrow d}
=
\mathbf I_{r_1}$
holds, and therefore
\begin{equation}
\mathbf Z_{r_1}(i)
=
\bar{\mathbf G}_{2\rightarrow d}(i),
\qquad
i=2,\dots,d.
\label{eq:Zr}
\end{equation}

The derivation for the remaining TT cores is analogous. For neighboring fibers associated with the left interpolation indices,
\[
X(\bar{\mathcal I}_{p-1}(i),:,\mathcal J_p)
=
X(\mathcal I_{p-1},:,\mathcal J_p)
\times_1
\mathbf Z_{l_p}(i),
\]
substituting the TT representation gives
\[
\bar{\mathbf G}_{1\rightarrow p-1}(i)\,
G_p
=
\mathbf Z_{l_p}(i)\,
\mathbf G_{1\rightarrow p-1}\,
G_p.
\]
Here, products involving \(G_p\) are understood after reshaping the \(p\)-th TT core into the matrix $\mathbf G_p\in\mathbb R^{r_{p-1}n_p\times r_p}.$
By construction, $\mathbf G_p(\mathcal I_p,:)=\mathbf I_{r_p},$
implying that \(\mathbf G_p\) has full column rank and therefore defines an injective linear map. Consequently, the common right factor \(\mathbf G_p\) can be eliminated. Since
$ \mathbf G_{1\rightarrow p-1}
=
\mathbf I_{r_{p-1}},$
it follows that
\begin{equation}
\mathbf Z_{l_p}(i)
=
\bar{\mathbf G}_{1\rightarrow p-1}(i),
\qquad
i=1,\dots,p-1.
\label{eq:Zl}
\end{equation}

Likewise, for neighboring fibers associated with the right interpolation indices,
\[
X(\mathcal I_{p-1},:,\bar{\mathcal J}_p(i))
=
X(\mathcal I_{p-1},:,\mathcal J_p)
\times_3
\mathbf Z_{r_p}(i)^{\top},
\]
the TT representation gives
\[
G_p\,
\bar{\mathbf G}_{p+1\rightarrow d}(i)
=
G_p\,
\mathbf G_{p+1\rightarrow d}\,
\mathbf Z_{r_p}(i),
\]
which, together with $\mathbf G_{p+1\rightarrow d}
=
\mathbf I_{r_p},$
yields
\begin{equation}
\mathbf Z_{r_p}(i)
=
\bar{\mathbf G}_{p+1\rightarrow d}(i),
\qquad
i=p+1,\dots,d.
\label{eq:Zrp}
\end{equation}

Finally, for the last TT core,
\[
X(\bar{\mathcal I}_{d-1}(i),:)
=
\mathbf Z_{l_d}(i)
X(\mathcal I_{d-1},:),
\]
the same argument gives
\begin{equation}
\mathbf Z_{l_d}(i)
=
\bar{\mathbf G}_{1\rightarrow d-1}(i),
\qquad
i=1,\dots,d-1.
\label{eq:Zld}
\end{equation}

The coefficient matrices associated with left dependencies are evaluated using the TT cores that have already been updated during the current forward sweep, whereas those associated with right dependencies are evaluated using the TT cores from the previous reverse sweep. This ordering follows naturally from the forward--backward structure of the \texttt{TT-CROSS-DEIM} algorithm and forms the basis of the fixed-point iteration described next.

\subsection{Fixed-point iteration for updating the interpolation coefficients}
\label{convergenceZ}

The interpolation coefficient matrices
$\mathbf Z_{r_1}(i)$,
$\mathbf Z_{l_p}(i)$,
$\mathbf Z_{r_p}(i)$,
and
$\mathbf Z_{l_d}(i)$
introduced in the previous subsection cannot be computed directly. Although
Eqs.~\eqref{eq:Zr}-\eqref{eq:Zld} provide explicit expressions for these matrices, they depend on
the neighboring TT contractions
$\bar{\mathbf G}$, which themselves are functions of the unknown TT
approximation at the new time step. Consequently, the TT cores and the
interpolation coefficient matrices must be determined simultaneously.

To this end, we employ a fixed-point iteration in which the TT cores and
the interpolation coefficient matrices are updated alternately. At the
beginning of each time step, the interpolation coefficient matrices are
initialized using the TT approximation from the previous time step.
Given the current estimate of the coefficient matrices, the closed
collocation systems
Eqs.~\eqref{Core1-3}, \eqref{Core2-3}, and \eqref{Core3-3}
are solved to obtain updated TT cores. These updated TT cores are then used
to evaluate the neighboring contractions
$\bar{\mathbf G}$, from which new interpolation coefficient matrices are
computed using
Eqs.~\eqref{eq:Zr}-\eqref{eq:Zld}. The procedure is repeated until
successive updates become self-consistent.

At convergence, the TT cores reproduce the neighboring fibers used to
construct the interpolation coefficient matrices, while the updated
interpolation coefficient matrices reproduce the same TT cores through the
closure relations
Eqs.~\eqref{eq:closure_first}--\eqref{eq:closure_last}. The resulting
fixed point therefore satisfies the original implicit collocation
equations.

To monitor convergence, we evaluate the residual of the implicit
collocation problem,
\begin{equation}
R
=
\mathcal A(\hat X^{k+1})
-
\mathcal B(\hat X^k),
\label{eq:residual_convergence}
\end{equation}
where $\hat X^{k+1}$ denotes the current TT approximation at the new
time step.

Because both $\hat X^{k+1}$ and $\hat X^k$ are represented in the TT
format, the residual can also be constructed directly in compressed form
without assembling any full-dimensional tensors. The relative residual is
then computed as
\begin{equation}
\label{eq:rel_res}
\varepsilon_R
=
\frac{\|R\|_F}
{\|\mathcal B(\hat X^k)\|_F},
\end{equation}
where all norms are evaluated using standard TT arithmetic.

The fixed-point iteration is terminated once
\begin{equation}
\varepsilon_R
<
\varepsilon_{\rm tol},
\end{equation}
for a prescribed tolerance $\varepsilon_{\rm tol}$. Throughout the
numerical examples considered in this work, the proposed iteration
converges rapidly, typically requiring fewer than 15 fixed-point
iterations per time step.

The dependency-elimination algorithm generates a sequence of intermediate \texttt{TT-CROSS-DEIM} approximations that remain on the fixed-rank TT manifold throughout the nonlinear solve. Consequently, every iteration corresponds to a valid low-rank tensor, and the implicit solution is obtained by converging directly on the manifold without tangent-space projections or optimization over TT cores. This is shown schematically in Figure \ref{fig:overview}.

\subsubsection{Rank adaptivity} \label{rankadaptivity}
The proposed time-integration scheme can be extended with a rank-adaptive mechanism that regulates the approximation error during the time advancement. For this purpose, we introduce the following error proxy:
\begin{equation} \label{ErrorCriterion}
\epsilon_q = \dfrac{\min(\boldsymbol \Sigma_q)}{\|\boldsymbol \Sigma_q \|_F}, \qquad q = 1,\dots,d-1,
\end{equation}
where $\boldsymbol \Sigma_q$ denotes the matrix of singular values obtained in the first step of Algorithm \ref{alg:CUR-DEIM} as $[\mathbf U,\boldsymbol \Sigma_q,\sim] = \texttt{SVD}(X(:,\mathcal J), r_q)$ for $q=1,\dots,d-1$, with $\mathcal J$ denoting the corresponding index set. The value of $\epsilon_q$ reflects the relative contribution of the $r_q$-th rank. The rank $r_p$ is either retained or modified to ensure that $\epsilon_q$ remains in a prescribed interval $\epsilon_l \le \epsilon_q \le \epsilon_u$, where $\epsilon_l$ and $\epsilon_u$ are user-defined lower and upper thresholds. If $\epsilon_q < \epsilon_l$, the rank is decreased by one, giving $r'_q = r_q - 1$, and the DEIM-based index set $\mathcal I_q$ is truncated to keep only the first $r'_q$ entries.
Conversely, if $\epsilon_q > \epsilon_u$, the rank is increased to $r'_q = r_q + 1$, which requires sampling additional indices and expanding $\mathcal I_q$.
Although rank modification could in principle be carried out during the iterations used to update the interpolation coefficient matrices $\mathbf Z$, we avoid doing so. Early iterations may produce inaccurate estimates of $\mathbf Z$, and adjusting the rank at that stage may degrade the overall accuracy of the time integration. Therefore, rank adaptivity is applied only once per time step, after the coefficient matrices have converged.

Since DEIM provides exactly as many sampling points as the number of columns of the left singular matrix $\mathbf U$, additional sampling is required to support rank increment. To obtain these extra points, we employ the GappyPOD+E procedure from \cite{BPStability}. Although other sparse sampling strategies could be used, GappyPOD+E has been shown to outperform methods such as random sampling or leverage-score sampling \cite{MMCUR}. Using GappyPOD+E, we append $h$ additional indices to the sets $\mathcal{I}_1,\dots,\mathcal{I}_{d-1}$, where $h\ge 1$ is an arbitrary number. Increasing $h$ generally reduces the approximation error; in practice, $h=5$ is sufficient to achieve the desired accuracy. Even though $h$ extra indices are added, the rank is increased by only one, since the SVD of the sampled function evaluations is truncated at $r+1$ (as shown in the first step of Algorithm \ref{alg:CUR-DEIM}).

\begin{algorithm}[t]
\fontsize{10pt}{11pt}\selectfont
\SetAlgoLined
\KwIn{Initial condition in TT form, initial interpolation coefficient matrices $\mathbf Z$} \vspace{1mm}
\KwOut{Solution of the TDE at the final time step in TT form} \vspace{1mm}
\For{k=1:max time steps}{\vspace{1mm}
    \While{ $\varepsilon_R<\varepsilon_{\rm tol}$}{\vspace{1mm}
        update $\mathbf Z_{r_1}(i)$ for $ i=2,\dots,d$ using Eq. (\ref{eq:Zr}) \\
        compute $X (:, \mathcal{J}_1)$ using Eq. (\ref{eq:sylv1}) \quad {\footnotesize \textit{\# use Newton iteration for nonlinear TDE}} \\
        compute first TT core: $G_1 \leftarrow X (:, \mathcal{J}_1)$  \vspace{1mm} \\ 
        \For{p=2:d-1}{
            update $\mathbf Z_{l_p}(i)$ for $i=1,\dots,p-1$ using Eq. (\ref{eq:Zl}) \\ 
            update $\mathbf Z_{r_p}(i)$ for $i=p+1,\dots,d$ using Eq. (\ref{eq:Zrp}) \\
            compute $X \bigl ( [\mathcal {I}_{p-1}, :] , \mathcal{J}_p \bigr )$ using Eq. (\ref{eq:sylv2}) \quad {\footnotesize \textit{\# use Newton iteration for nonlinear TDE}}\\
            compute $p$-th TT core: $G_p \leftarrow X \bigl ( [\mathcal {I}_{p-1}, :] , \mathcal{J}_p \bigr )$   \\
            } \vspace{1mm}
        update $\mathbf Z_{l_d}(i)$ for $ i=1,\dots,d-1$ using Eq. (\ref{eq:Zld}) \\
        compute $X (\mathcal{I}_{d-1}, :)$ using Eq. (\ref{eq:sylv3}) \quad {\footnotesize \textit{\# use Newton iteration for nonlinear TDE}}\\
        compute last TT core: $G_d \leftarrow X (\mathcal{I}_{d-1}, :)$   \\ 
        \vspace{1.5 mm}
        compute the relative residual $\varepsilon_R$ using Eq. (\ref{eq:rel_res})\\
        \vspace{0.5 mm}
        reassign the nested indices for next reverse sweep: $\mathcal{J}_1 \leftarrow \mathcal{I}_{d-1};~ \mathcal{J}_2 \leftarrow \mathcal{I}_{d-2}; ~\dots ~; ~ \mathcal{J}_{d-1} \leftarrow \mathcal{I}_1$ \hspace{5mm}  \\
        \vspace{0.5 mm}
        permute TT cores' indices and reassign them: $G_1 \leftarrow \texttt{permute} \bigl (G_d, [3,2,1] \bigr); \quad \dots \quad; G_d \leftarrow \texttt{permute} \bigl (G_1, [3,2,1] \bigr)$
        }
    }
    \caption{Implicit time integration of TDEs using TT-CROSS-DEIM}
\label{alg:Implicit-TT-CROSS-DEIM}
\end{algorithm}

\subsubsection{Extension to nonlinear tensor differential equations}
\label{ExtendNonlinear}

The dependency-elimination framework developed above extends naturally to tensor differential equations with general nonlinearities. Once the neighboring-fiber dependencies have been eliminated, the implicit collocation equations for the \texttt{TT-CROSS-DEIM} subtensors remain closed, but are generally nonlinear. Consequently, each TT subtensor is obtained by solving a low-dimensional nonlinear system.

Let
\[
\mathcal F(X)=0
\]
denote the nonlinear system associated with one of the TT subtensors (e.g.,
\(X(:,\mathcal J_1)\),
\(X(\mathcal I_{p-1},:,\mathcal J_p)\), or
\(X(\mathcal I_{d-1},:)\)).
We solve this system using Newton's method,
\begin{equation}
\mathbf J(X^{(m)})\,\delta X
=
-\mathcal F(X^{(m)}),
\qquad
X^{(m+1)}
=
X^{(m)}
+
\delta X,
\end{equation}
where \(\mathbf J\) denotes the Jacobian of \(\mathcal F\).

An important advantage of the proposed framework is that the nonlinear operator is never assembled in full tensor form. Instead, the residual and, when needed, Jacobian-vector products are evaluated directly from pointwise evaluations of the nonlinear operator on the DEIM-selected fibers. Consequently, arbitrary nonlinearities—including polynomial, exponential, logarithmic, fractional, and other nonlinear constitutive laws—can be incorporated without modifying the underlying TT algorithm.

Because each Newton solve involves only the DEIM-selected subtensor, the resulting nonlinear systems are relatively small. They may therefore be solved using either conventional Newton's method with an explicitly assembled Jacobian or Jacobian-free Newton--Krylov (JFNK) methods for larger problems \cite{Knoll2004}. Within each outer iteration that updates the interpolation coefficient matrices $\mathbf Z$, Newton iterations are performed until the nonlinear residual satisfies a prescribed tolerance.

In all numerical examples considered in this work, Newton's method converged rapidly, typically requiring only two to five iterations. Algorithm \ref{alg:Implicit-TT-CROSS-DEIM} summarizes the main steps for implicit time integration of TDEs on the low-rank TT manifolds.

\subsubsection{Computational Complexity} \label{ComputationalComplexity}
In this section, we analyze the computational complexity of the main components of the proposed algorithm for the implicit time integration of TDEs.  The following complexity analysis is presented for the class of implicit TDEs considered in this work. In particular, we assume that the underlying TDE can be expressed in the form of Eq.~(\ref{eq:linearTDE}), for which the local TT subproblems derived in Sections \ref{sec:interp_coef}--\ref{sec:interp_coef} apply. The resulting estimates quantify the computational complexity of the proposed implicit TT-cross framework, excluding the problem-dependent cost of evaluating the underlying differential operator. 

The overall computational cost is obtained by examining each stage of the algorithm individually.
\begin{itemize}
    \item \textbf{Updating the interpolation matrices:} Assuming that the one-dimensional operators have bounded sparsity, so that the neighboring index-set sizes satisfy $\bar r_p(i)=\mathcal{O}(r)$, the computation of a single interpolation coefficient matrix requires $\mathcal{O}(dr^3)$ operations. Specifically, the interpolation coefficient matrix is assembled through a sequence of contractions of the TT cores evaluated at the DEIM-selected points. Each contraction of the accumulated product with the next TT core incurs a computational cost of $\mathcal{O}(r^3)$ over all $r$ DEIM points. Since constructing one interpolation coefficient matrix requires $\mathcal{O}(d)$ such contractions, its total computational cost is $\mathcal{O}(dr^3)$. During each fixed-point iteration, interpolation coefficient matrices are computed for $\mathcal{O}(d)$ operator terms associated with each of the $d$ TT cores. Consequently, the total cost of updating all interpolation coefficient matrices during one iteration is $\mathcal{O}(d^3r^3)$.
    

    \item \textbf{Evaluation of the right-hand sides:}  Each evaluation of the right-hand sides of Eqs. (\ref{eq:sylv1}, \ref{eq:sylv2}, \ref{eq:sylv3}), namely $\mathbf{F}_1$, $\mathbf{F}_p$, and $\mathbf{F}_d$, requires $\mathcal{O}(dr^3 + nr^2)$ operations. Consequently, evaluating the right-hand sides for all $d$ TT cores during one iteration results in a total computational cost of $\mathcal{O}(d^2r^3+dnr^2)$.

    \item \textbf{Singular value decompositions:}  The most computationally expensive SVD is performed for the middle TT cores, where its complexity is $\mathcal{O}(nr^3)$. Since an SVD is computed for each TT core, the total computational cost of all SVD computations during one iteration is $\mathcal{O}(dnr^3)$.

   \item \textbf{Solving the local subproblems:}
The computational cost of solving the local subproblems depends on both the
underlying TDE and the numerical solver employed. When a direct dense linear
solver is used, solving all $d$ local subproblems requires a total cost of
$\mathcal{O}(dn^{3}r^{6})$. For linear TDEs, the local subproblems can instead
be formulated as Sylvester equations and solved using the Bartels--Stewart
algorithm with a total computational cost of
$\mathcal{O}(dn^{3}r^{3})$.
For larger values of $n$ and $r$, iterative Krylov subspace methods such as
GMRES can significantly reduce the computational cost. This is particularly
advantageous for nonlinear TDEs, where the local subproblems cannot be
reformulated as Sylvester equations and the use of direct linear solvers
becomes prohibitively expensive. Assuming that $k$ GMRES iterations are
required for each local solve, each Jacobian--vector product requires
$\mathcal{O}(n^{2}r^{2}+nr^{3})$ operations, corresponding to the left and
right matrix multiplications in the local operator. Consequently, solving all
$d$ local subproblems requires a total computational cost of $\mathcal{O}\!\left(kd\left(n^{2}r^{2}+nr^{3}\right)\right)$.
\end{itemize}

From both the above complexity estimates and our numerical experiments, solving the local subproblems is the dominant computational cost of the proposed algorithm, while the costs associated with updating the interpolation matrices, evaluating the right-hand sides, and computing the SVDs are comparatively small. Consequently, the overall performance is largely determined by the efficiency of the local linear or nonlinear solver.

\section{Demonstrations} \label{sec:Demo}
\subsection{3D heat equation} \label{sec:3Dheat}

As an initial test case, we consider the three-dimensional heat equation with a source term,
\begin{equation}
    \label{eq:3DheatEq}
    \dfrac{\partial v(\mathbf x,t)}{\partial t}
    =
    \alpha \nabla^2v(\mathbf x,t)
    +
    q(\mathbf x,t),
    \qquad
    \mathbf x\in[0,1],
\end{equation}
where $\alpha$ denotes the diffusion coefficient and $q(\mathbf x,t)$ is the heat source. This example is designed to examine four aspects of the proposed methodology: (i) the accuracy of the low-rank implicit solution, (ii) the convergence behavior of the dependency-elimination iteration, (iii) preservation of the temporal order of the underlying implicit time integrator, and (iv) the dependence of the fixed-point iteration on the spatial resolution.

The problem is equipped with homogeneous Dirichlet boundary conditions, and the initial condition is prescribed as
\begin{equation}
    v(\mathbf x,0)
    =
    \sum_{i=1}^{6}
    \theta_i
    \exp
    \Big(
    -60
    \big(
    (x_1-c_i)^2
    +
    (x_2-c_i)^2
    +
    (x_3-c_i)^2
    \big)
    \Big),
\end{equation}
where
$\theta_1=1$,
$\theta_2=0.5$,
$\theta_3=0.3$,
$\theta_4=0.25$,
$\theta_5=0.2$,
$\theta_6=10^{-4}$,
and
$c_1=0.5$,
$c_2=0.55$,
$c_3=0.45$,
$c_4=0.6$,
$c_5=0.4$,
$c_6=0.65$.
The spatial domain is discretized using second-order central finite differences with
$n_1=n_2=n_3=60$ grid points, and the diffusion coefficient is set to $\alpha=1$. The heat source is chosen as
\[
q(\mathbf x,t)=1-0.5\sin(2\pi t),
\]
which is spatially uniform. Time integration is performed over
$t\in[0,1]$
using the implicit Euler method with
$\Delta t=0.01$.
The resulting tensor differential equation is solved on the TT manifold using fixed TT ranks
$r=5$
and
$r=8$,
where
$r=r_1=r_2=r_3$.

Accuracy is measured using the relative error
\begin{equation}
\label{eq:err_rel}
\mathcal E(t)
=
\frac{\|
\hat V(t)-V(t)
\|_F}
{\|
V(t)
\|_F},
\end{equation}
where
$\hat V(t)$
denotes the TT approximation and
$V(t)$
is the full-order solution computed with the same spatial and temporal discretizations.

\begin{figure}[!t]
     \centering
        \includegraphics[width=0.8\linewidth]{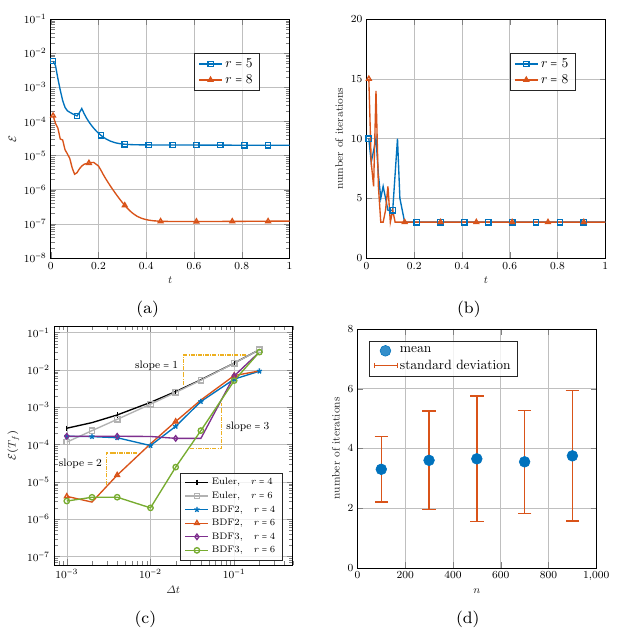} \vspace{-3mm}
     \caption{Three-dimensional heat equation (case 1); (a) relative error evolution with two fixed ranks; (b) number of the required iterations for convergence of the $\mathbf Z$ matrices over time associated with the plots of Figure \ref{fig:3DHeat1}a; (c) error at the final time ($\mathcal{E}(T_f)$) versus time step size ($\Delta t$). The implicit Euler method, BDF2, and BDF3 are used for time integration with various ranks $r=4$, $r=6$; (d) mean and standard deviation of the number of iterations for convergence of the coefficient matrices versus $n$.}
     \label{fig:3DHeat1}
\end{figure}

Figure~\ref{fig:3DHeat1}a shows the evolution of the relative error for the two TT ranks. In both cases the error decreases rapidly during the initial transient before reaching a nearly constant level determined primarily by the prescribed TT rank. As expected, increasing the rank significantly reduces the approximation error. 

The convergence of the dependency-elimination iteration is illustrated in Figure~\ref{fig:3DHeat1}b. During the first few time steps the coefficient matrices require between approximately 10 and 15 fixed-point iterations to converge. As the transient decays, however, the required number of iterations rapidly decreases and subsequently stabilizes at approximately three iterations per time step. This behavior reflects the temporal continuity of the solution, whereby the converged interpolation coefficients from one time step provide an increasingly accurate initial guess for the next. The figure therefore validates the efficiency of the proposed fixed-point iteration, which constitutes the central algorithmic component of the present work.

A key advantage of performing the low-rank approximation after temporal discretization is that the proposed formulation preserves the temporal accuracy of the underlying implicit integrator. Figure~\ref{fig:3DHeat1}c reports the error at the final time,
$\mathcal E(T_f)$,
as a function of the time-step size for the implicit Euler, BDF2, and BDF3 schemes using ranks
$r=4$
and
$r=6$.
The observed convergence rates coincide with the theoretical first-, second-, and third-order accuracy of the corresponding full-order discretizations until the error reaches the low-rank approximation floor. In particular, the higher-order BDF3 scheme reaches this saturation level using considerably larger time steps than BDF2. These results demonstrate that the proposed dependency-elimination procedure does not degrade the temporal accuracy of the underlying implicit discretization.

Finally, we investigate how the convergence of the interpolation-coefficient iteration depends on the spatial resolution. Since the size of each coefficient matrix $\mathbf Z$ increases with the number of grid points, one might expect the fixed-point iteration to become increasingly difficult for larger problems. To examine this effect, we repeat the simulation while varying the one-dimensional grid size from $n=100$ to $n=900$. Figure~\ref{fig:3DHeat1}d reports the mean and standard deviation of the number of fixed-point iterations. Remarkably, both quantities remain nearly unchanged over the entire range of spatial resolutions. Although increasing $n$ substantially enlarges the interpolation coefficient matrices, no noticeable deterioration in convergence is observed. This suggests that the convergence behavior of the proposed dependency-elimination iteration is governed primarily by the underlying low-rank structure rather than by the ambient tensor dimension.

Although the proposed algorithm is designed to handle general nonlinear TDEs, it is also applicable to linear problems. Therefore, it is instructive to compare its performance against established state-of-the-art methods for the implicit solution of linear TDEs. In this section, we consider two widely used approaches: ALS \cite{Zach2012, Rohwedder2013} and AMEn \cite{Dolgov2014}. Since two-site DMRG is algorithmically equivalent to two-site ALS for solving linear systems in the TT format, we expect the computational performance of the proposed method relative to DMRG to be similar to that observed for ALS. These methods are well-established techniques for solving linear systems in the TT format and provide a meaningful benchmark for assessing the performance of the proposed algorithm.

The one-site ALS algorithm operates with a prescribed TT rank that remains fixed throughout the iterations. Consequently, it cannot adapt the tensor ranks to reflect the evolving complexity of the solution during the time integration. This limitation is addressed by the two-site ALS algorithm, in which two neighboring TT cores are merged and optimized simultaneously. The merged core is subsequently factorized into two TT cores using a singular value decomposition, allowing the intermediate TT rank to increase or decrease automatically according to a prescribed truncation tolerance. While this procedure provides rank adaptivity, the repeated merging and SVD operations substantially increase the computational cost compared with one-site ALS. The AMEn algorithm was proposed to alleviate this drawback by introducing rank adaptivity without merging neighboring TT cores. Instead, it enriches the local approximation space using information extracted from the residual, thereby achieving adaptive rank growth with a lower computational overhead.

\begin{figure}[!t]
     \centering
    \includegraphics[width=0.8\linewidth]{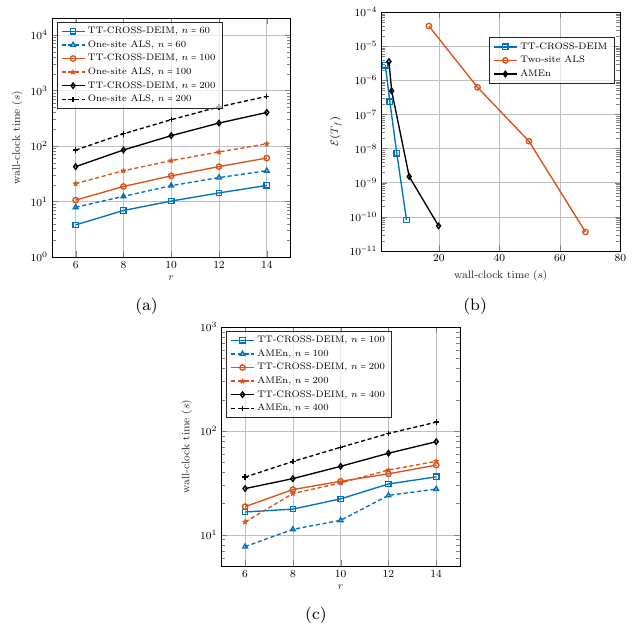} \vspace{-3mm}
     \caption{Three-dimensional heat equation (case 1); (a) CPU time versus rank for \texttt{TT-CROSS-DEIM} and ALS under different choices of $n$; (b) relative error at the final time step versus CPU time for rank adaptive versions of \texttt{TT-CROSS-DEIM}, two-site ALS,  and AMEn; (c) CPU time versus rank for \texttt{TT-CROSS-DEIM} and AMEn under different choices of $n$ and increased dimension of the problem ($d=30$).}
     \label{fig:3DHeat1-compare}
\end{figure}

We first compare the proposed \texttt{TT-CROSS-DEIM} algorithm with the one-site ALS method under fixed TT ranks. Figure \ref{fig:3DHeat1-compare}a reports the wall-clock time required to perform 30 iterations within a single implicit time step for a fixed time step size $\Delta t$. The comparison is carried out for three spatial discretizations, namely $n=60$, $n=100$, and $n=200$, over a range of prescribed TT ranks. As shown in Figure \ref{fig:3DHeat1-compare}a, the computational cost of both methods increases with the TT rank and becomes larger as the spatial resolution increases. However, \texttt{TT-CROSS-DEIM} consistently outperforms one-site ALS across all tested configurations, requiring approximately half the execution time for every combination of $r$ and $n$.

Figure \ref{fig:3DHeat1-compare}b compares the cost-error trade-offs of the rank-adaptive versions of \texttt{TT-CROSS-DEIM}, two-site ALS, and AMEn. The results indicate that \texttt{TT-CROSS-DEIM} and AMEn exhibit comparable accuracy for a given computational cost, with both methods rapidly reducing the relative error to near machine precision. In contrast, two-site ALS requires substantially more wall-clock time to achieve the same level of accuracy. This additional cost stems primarily from the repeated core-merging and SVD operations required to maintain rank adaptivity. Overall, these results demonstrate that \texttt{TT-CROSS-DEIM} combines the computational efficiency of fixed-rank algorithms with the adaptive capabilities of modern rank-adaptive methods, yielding an excellent balance between accuracy and computational cost. Moreover, unlike ALS and AMEn, which are limited to linear TDEs, \texttt{TT-CROSS-DEIM} naturally extends to nonlinear TDEs, making it applicable to a substantially broader class of problems.

To assess the scalability of AMEn and \texttt{TT-CROSS-DEIM}, we increase the dimensionality of the current problem from three to thirty while keeping all other problem parameters unchanged. Unlike the experiments reported in Figures \ref{fig:3DHeat1-compare}a and \ref{fig:3DHeat1-compare}b, where the local subproblems were solved using a direct Sylvester solver, the local linear systems in this experiment are solved iteratively using the GMRES algorithm. Figure \ref{fig:3DHeat1-compare}c presents the wall-clock time required to perform 30 iterations within a single implicit time step for a fixed time step size $\Delta t$. The computational performance is evaluated for three spatial discretizations, with $n=100$, $n=200$, and $n=400$, while varying the prescribed TT rank. As shown in Figure \ref{fig:3DHeat1-compare}c, the relative computational performance of the two methods depends on both the spatial resolution and the prescribed TT rank. For the coarsest discretization ($n=100$), AMEn consistently requires less computational time than \texttt{TT-CROSS-DEIM} for all tested ranks. When the discretization is refined to $n=200$, the performance gap narrows, and the two methods exhibit nearly identical computational costs for moderate ranks. As the prescribed TT rank increases, however, \texttt{TT-CROSS-DEIM} becomes slightly more efficient than AMEn. For the finest discretization ($n=400$), \texttt{TT-CROSS-DEIM} consistently outperforms AMEn across the entire range of tested ranks, and the performance gap widens as the TT rank increases.  The improved performance of \texttt{TT-CROSS-DEIM} for larger problems can be explained by the dominant computational cost of both algorithms, namely the solution of the local subproblems. As the spatial resolution and TT rank increase, the cost of these local solves grows for both methods. However, unlike AMEn, the proposed method does not require residual-based basis enrichment or the associated orthogonalization and subspace updates. Since the remaining operations, such as updating the interpolation matrices, represent only a small fraction of the overall runtime, the computational advantage of \texttt{TT-CROSS-DEIM} becomes increasingly pronounced for larger problems.

For the second case, all parameters remain unchanged except for the heat source,
\[
q(\mathbf x,t)
=
\exp
\Big(
-80
\big(
(x_1-0.3)^2
+
(x_2-0.2)^2
+
(x_3-0.7)^2
\big)
\Big),
\]
which is time independent. The simulation interval is extended to
$t\in[0,100]$,
and the implicit Euler method is employed with a substantially larger time-step size,
$\Delta t=10$.
The problem is solved using TT ranks
$r=7$
and
$r=10$.

\begin{figure}[htbp]
     \centering
     \includegraphics[width=0.8\linewidth]{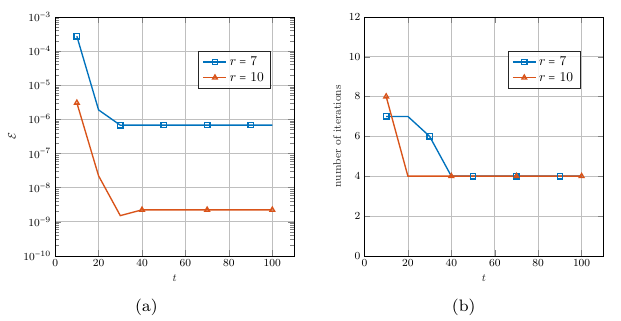} \vspace{-3mm}
     \caption{Three-dimensional heat equation (case 2); (a) relative error evolution with two fixed ranks; (b) number of the required iterations for convergence of the $\mathbf Z$ matrices over time associated with the plots of Figure \ref{fig:3DHeat2}a.}
     \label{fig:3DHeat2}
\end{figure}

The relative errors shown in Figure~\ref{fig:3DHeat2}a again decrease as the TT rank increases and remain essentially constant over the entire simulation interval. More importantly, the corresponding convergence histories in Figure~\ref{fig:3DHeat2}b indicate that only four to eight fixed-point iterations are required despite the time step being increased by three orders of magnitude relative to the first test case. Neither the approximation error nor the convergence behavior deteriorates as the time step becomes much larger. These results demonstrate that the proposed dependency-elimination framework remains robust for large implicit time steps and is therefore well suited for the efficient computation of long-time and steady-state solutions.

\subsection{Thirty-dimensional heat equation}
\label{sec:nD-heat}

As the second demonstration, we consider the 30-dimensional heat equation
\begin{equation}
\label{eq:nDheatEq}
\dfrac{\partial v(\mathbf x,t)}{\partial t}
=
\alpha\nabla^2v(\mathbf x,t)
+
q(\mathbf x,t),
\qquad
\mathbf x\in[0,5],
\end{equation}
with diffusion coefficient $\alpha=0.5$. This example is intended to demonstrate the capability of the proposed method for solving tensor differential equations whose ambient dimension is far beyond the reach of conventional discretization techniques. In particular, we investigate the convergence of the rank-adaptive algorithm and the evolution of the TT ranks throughout the simulation.

Zero Dirichlet boundary conditions are imposed, and the initial condition is chosen as
\begin{equation}
v(\mathbf x,0)
=
\exp
\Big(
-2
\sum_{i=1}^{30}(x_i-3.25)^2
\Big)
+
\exp
\Big(
-2
\sum_{i=1}^{30}(x_i-1.75)^2
\Big),
\end{equation}
while the heat source is given by
\[
q(\mathbf x,t)
=
(2-e^{-t})
\exp
\Big(
-2
\sum_{i=1}^{30}(x_i-2.5)^2
\Big).
\]

Each spatial coordinate is discretized using second-order central finite differences with
\[
n_1=n_2=\cdots=n_{30}=70,
\]
and implicit Euler time integration is performed over
$t\in[0,1]$
using
$\Delta t=0.05$. The corresponding full-order discretization contains
\[
70^{30}
\approx
2.3\times10^{55}
\]
degrees of freedom, making direct numerical simulation completely infeasible. Consequently, no full-order reference solution is available for assessing the accuracy of the low-rank approximation. Instead, we evaluate convergence by examining low-dimensional marginal functions of the solution. Specifically, we consider the two-dimensional marginal
\begin{equation}
\label{eq:marginal}
p(x_9,x_{10},t)
=
\int_{[0,5]^{28}}
v(\mathbf x,t)
\,d\mathbf x_1
\cdots
d\mathbf x_8
d\mathbf x_{11}
\cdots
d\mathbf x_{30},
\end{equation}
which can be computed efficiently from the TT representation.

Rank adaptivity is employed throughout the simulation. Four pairs of adaptivity thresholds,
$
(\epsilon_l,\epsilon_u),
$
denoted by
$\epsilon_1$,
$\epsilon_2$,
$\epsilon_3$,
and
$\epsilon_4$,
are considered:
\[
\epsilon_1:
(10^{-8},10^{-7}),
\quad
\epsilon_2:
(10^{-9},10^{-8}),
\quad
\epsilon_3:
(10^{-10},10^{-9}),
\quad
\epsilon_4:
(10^{-11},10^{-10}).
\]
Reducing these thresholds permits additional TT ranks to be retained, thereby increasing the approximation accuracy.

Since no full-order solution is available, the solution obtained with the tightest tolerances,
$\epsilon_4$,
is adopted as the reference solution. For each remaining threshold pair, we compute the difference between the corresponding marginal function and that obtained with
$\epsilon_4$,
normalized by the range of the reference marginal,
\[
\max(p_{\epsilon_4})-\min(p_{\epsilon_4}).
\]
\begin{figure}[!t]
    \centering
    \includegraphics[scale=1]{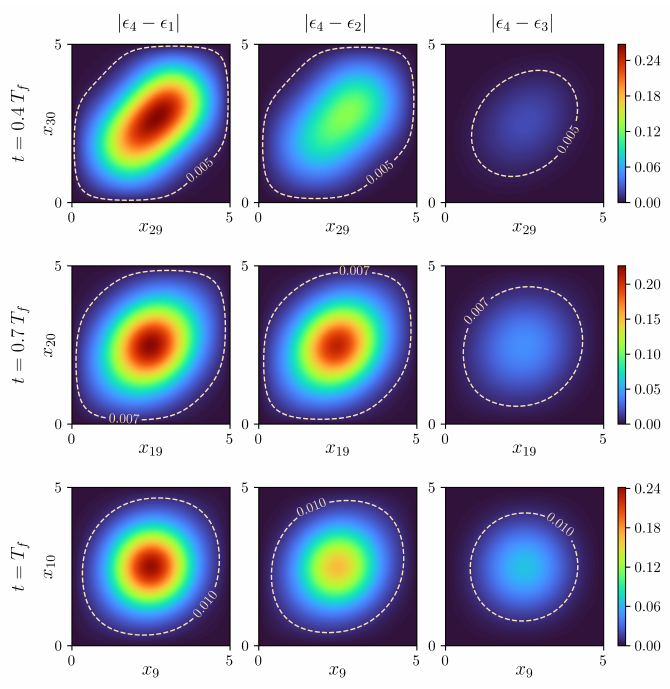} \vspace{-3mm}
     \caption{30D heat equation: relative difference of the marginal plots of joint positions of $(\bm x_{29}, \bm x_{30})$ at $0.4\:T_f$ (first row), $(\bm x_{19}, \bm x_{20})$ at $t=0.7\:T_f$ (second row), and $(\bm x_9, \bm x_{10})$ at $t=T_f$ (third row); obtained by solving the TDE with $\epsilon_{4}$ and $\epsilon_{1}$ (Left column); $\epsilon_{4}$ and $\epsilon_{2}$ (Middle column); $\epsilon_{4}$ and $\epsilon_{3}$ (Right column).}
     \label{fig:nDheat_convergance}
\end{figure}

Figure~\ref{fig:nDheat_convergance} presents the resulting convergence study. The three columns report the normalized differences between the marginal functions obtained using
$\epsilon_1$,
$\epsilon_2$,
and
$\epsilon_3$
relative to
$\epsilon_4$,
respectively, evaluated for the coordinate pairs $(\bm x_{29},\bm x_{30})$ at $t=0.4T_f$, $(\bm x_{19},\bm x_{20})$ at $t=0.7T_f$, and $(\bm x_{9},\bm x_{10})$ at $t=T_f$.
A clear monotonic reduction in the error is observed as the adaptivity tolerances are tightened, demonstrating convergence of the rank-adaptive algorithm despite the absence of a full-order reference solution. This experiment therefore provides strong evidence that the proposed implicit TT-cross formulation converges systematically as additional low-rank degrees of freedom are introduced.

\begin{figure}
     \centering
     \includegraphics[width=0.9\linewidth]{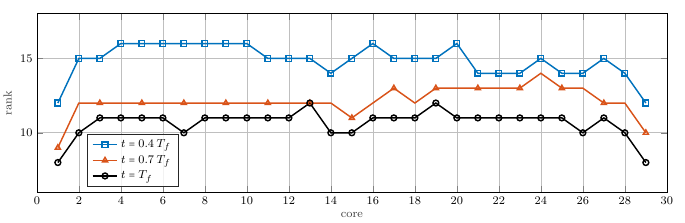} \vspace{-3mm}
     \caption{30D heat equation: all tensor train ranks at different time steps obtained by solving the system with $\epsilon_4: \:\epsilon_{u} = 10^{-10}$ and $\epsilon_{l} = 10^{-11}$.}
     \label{fig:ranks_nDheat}
\end{figure}

An additional advantage of the adaptive formulation is illustrated in Figure~\ref{fig:ranks_nDheat}, which reports the TT ranks at
$t=0.4T_f$,
$t=0.7T_f$,
and
$t=T_f$
for the threshold pair
$\epsilon_4$.
Since a 30-dimensional tensor possesses 29 TT ranks, the figure provides a detailed view of how the complexity of the solution evolves over time.
The ranks remain nonuniform across the tensor network and gradually decrease as the simulation progresses. Physically, this behavior reflects the diffusive nature of the heat equation: as the transient components decay, the solution becomes progressively smoother and therefore admits a lower-rank representation.
The ability of the algorithm to automatically reduce the TT ranks during the simulation illustrates one of the principal advantages of combining the proposed implicit integrator with rank adaptivity, namely that the computational complexity naturally tracks the intrinsic complexity of the evolving solution rather than remaining fixed throughout the computation.

\subsection{3D nonlinear advection--reaction equation}
\label{sec:3D-AR}

We next consider a three-dimensional nonlinear advection--reaction equation,
\begin{equation}
    \label{eq:3DAR}
    \dfrac{\partial v(\mathbf x,t)}{\partial t}
    =
    \alpha \nabla v(\mathbf x,t)
    +
    \beta
    \dfrac{v(\mathbf x,t)}{1+v(\mathbf x,t)^2},
    \qquad
    \mathbf x\in[-5,5],
\end{equation}
where $\alpha=-0.5$ is the advection velocity and $\beta=-0.5$ controls the nonlinear reaction term. This example is designed to assess the performance of the proposed method for a nonlinear TDE requiring Newton iterations inside the dependency-elimination loop. The reaction term
is deliberately chosen to be non-polynomial. Unlike polynomial nonlinearities, such nonlinear operators cannot be efficiently represented through preassembled multilinear tensor operators, whereas the \texttt{TT-CROSS-DEIM} formulation requires only pointwise evaluations of the nonlinear term and therefore naturally accommodates general nonlinear functions.

Periodic boundary conditions are imposed, and the initial condition is
\[
v(\mathbf x,0)
=
\exp\!\left(
-\bigl(x_1^2+x_2^2+x_3^2\bigr)
\right)
+
\exp\!\left(
-\bigl((x_1+2)^2+(x_2+2)^2+(x_3+2)^2\bigr)
\right).
\]
The domain is discretized using second-order central finite differences with
$n_1=n_2=n_3=70$ grid points. Time integration is performed with implicit Euler over
$t\in[0,5]$ using $\Delta t=0.01$.

\begin{figure}
     \centering
    \includegraphics[width=0.8\linewidth]{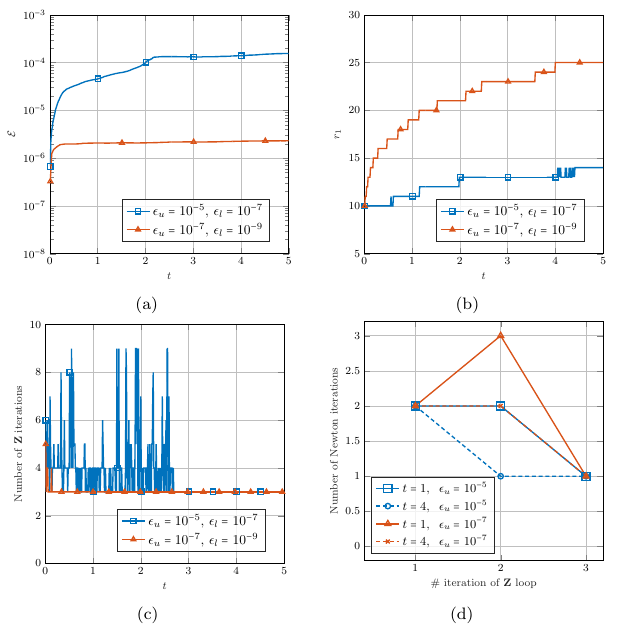} \vspace{-3mm}

     \caption{Three-dimensional advection reaction equation; (a) relative error evolution with two rank-adaptivity thresholds; (b) evolution of $r_1$ associated with the rank-adaptive systems plotted in Figure \ref{fig:3DAR}a; (c) number of the required iterations for convergence of the $\mathbf Z$ loops over time associated with the plots of Figure \ref{fig:3DAR}a; (d) number of the required iterations for convergence of the Newton iterations of the middle core ($G_2$) at two specific time steps ($t=1,4$). The horizontal axis shows the first, second, and third $\mathbf Z$ loops at the mentioned time-steps.}
     \label{fig:3DAR}
\end{figure}

The FOM is used as the reference solution. The nonlinear TDE is solved using the proposed implicit \texttt{TT-CROSS-DEIM} method with rank adaptivity. Two sets of adaptivity thresholds are considered:
\[
(\epsilon_u,\epsilon_l)=(10^{-5},10^{-7}),
\qquad
(\epsilon_u,\epsilon_l)=(10^{-7},10^{-9}).
\]
Figure~\ref{fig:3DAR}a shows the relative error between the low-rank solution and the FOM. As expected, tightening the adaptivity thresholds retains higher TT ranks and reduces the approximation error. The corresponding evolution of the first TT rank, shown in Figure~\ref{fig:3DAR}b, illustrates how the adaptive procedure adjusts the representation during the simulation.

Figure~\ref{fig:3DAR}c reports the number of fixed-point iterations required for convergence of the interpolation coefficient matrices $\mathbf Z$. For $(\epsilon_u,\epsilon_l)=(10^{-5},10^{-7})$, the iteration count ranges from 3 to 9, while for $(\epsilon_u,\epsilon_l)=(10^{-7},10^{-9})$, it ranges from 3 to 5. Thus, even though the governing equation contains a nonlinear reaction term, the dependency-elimination iteration remains rapidly convergent. This indicates that the closure iteration is not substantially degraded by the nonlinear reaction term in this example.

For nonlinear TDEs, each fixed-point iteration requires Newton solves for the \texttt{TT-CROSS-DEIM} subtensors, as described in Section~\ref{ExtendNonlinear}. To examine the cost of these inner nonlinear solves, we select two representative time instances, $t=1$ and $t=4$. At both times, three fixed-point iterations are sufficient for convergence of the coefficient matrices for both adaptivity thresholds. Figure~\ref{fig:3DAR}d shows the number of Newton iterations required to compute the middle TT core $G_2$ within each of these fixed-point iterations. Across all cases, only a small number of Newton iterations is required. These results demonstrate that the proposed framework remains efficient for nonlinear implicit problems and the outer dependency-elimination iteration converges rapidly.

\subsection{Twenty-dimensional nonlinear advection--reaction equation}
\label{sec:nD-AR}

As the final demonstration, we consider a twenty-dimensional nonlinear advection--reaction equation,
\begin{equation}
    \label{eq:nDAR}
    \dfrac{\partial v(\mathbf x,t)}{\partial t}
    =
    \alpha \nabla v(\mathbf x,t)
    +
    \beta
    \frac{v(\mathbf x,t)}
    {1+v(\mathbf x,t)^2},
    \qquad
    \mathbf x\in[-5,5],
\end{equation}
where $\alpha=-0.5$ is the advection velocity and $\beta=-0.5$ is the coefficient of the nonlinear reaction term. This example combines the two principal challenges addressed by the proposed methodology. First, the advection operator introduces non-diagonal couplings between neighboring fibers, making the dependency-elimination framework essential for obtaining a closed implicit collocation system. Second, the nonlinear reaction term is deliberately chosen to be non-polynomial.

Periodic boundary conditions are imposed, and the initial condition is
\[
v(\mathbf x,0)
=
\exp\!\Big(
-\big((x_1-0.5)^2+\cdots+(x_{20}-0.5)^2\big)
\Big)
+
\exp\!\Big(
-\big((x_1+1.5)^2+\cdots+(x_{20}+1.5)^2\big)
\Big).
\]
Each spatial dimension is discretized using second-order central finite differences with
$n_1=\cdots=n_{20}=100$ grid points. Time integration is performed using the implicit Euler method over the interval
$t\in[0,2]$
with a uniform time step
$\Delta t=0.01$.

The corresponding full-order discretization contains
\[
100^{20}=10^{40}
\]
degrees of freedom, rendering a direct numerical simulation entirely infeasible. Consequently, no full-order reference solution is available. Following the previous high-dimensional example, convergence is assessed through the time-dependent marginal functions defined in Eq.~\eqref{eq:marginal}.

To investigate convergence with respect to the adaptive TT representation, four pairs of adaptivity thresholds,
$
(\epsilon_l,\epsilon_u),
$
denoted by
$\epsilon_1$,
$\epsilon_2$,
$\epsilon_3$,
and
$\epsilon_4$,
are considered:
\[
\epsilon_1:
(10^{-3},10^{-2}),
\qquad
\epsilon_2:
(10^{-4},10^{-3}),
\qquad
\epsilon_3:
(10^{-5},10^{-4}),
\qquad
\epsilon_4:
(10^{-6},10^{-5}).
\]
The solution computed using $\epsilon_4$ is taken as the reference solution. For each of the remaining threshold pairs, the marginal functions are compared with those obtained using $\epsilon_4$, and the differences are normalized by the range of the corresponding reference marginal function, $\max(p_{\epsilon_4})-\min(p_{\epsilon_4})$.

\begin{figure}[!t]
    \centering
    \includegraphics[scale=1]{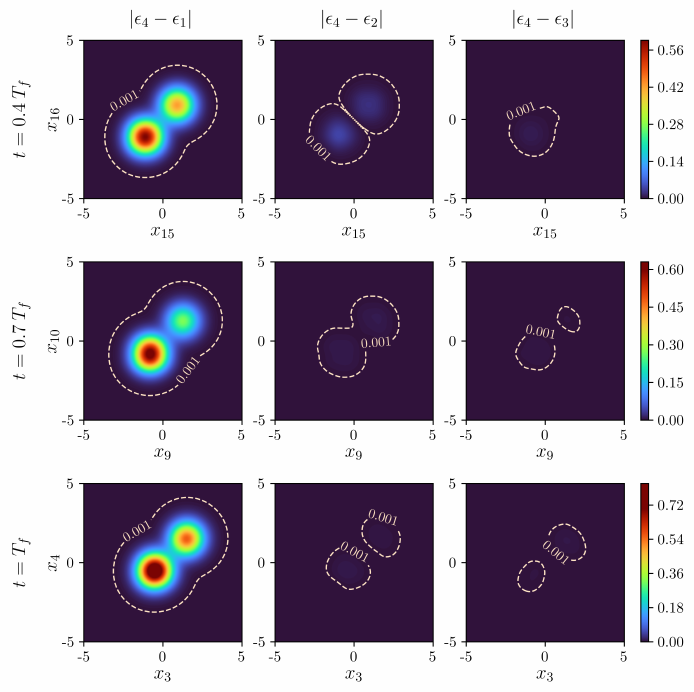} \vspace{-3mm}
     \caption{20D advection reaction equation: relative difference of the marginal plots of joint positions of $(\bm x_{15}, \bm x_{16})$ at $0.4\:T_f$ (first row), $(\bm x_{9}, \bm x_{10})$ at $t=0.7\:T_f$ (second row), and $(\bm x_3, \bm x_{4})$ at $t=T_f$ (third row); obtained by solving the TDE with $\epsilon_{4}$ and $\epsilon_{1}$ (Left column);  $\epsilon_{4}$ and $\epsilon_{2}$ (Middle column); $\epsilon_{4}$ and $\epsilon_{3}$ (Right column).}
     \label{fig:nDAR_convergance}
\end{figure}

Figure~\ref{fig:nDAR_convergance} summarizes the convergence study. The left, middle, and right columns display the normalized differences between the marginal functions computed with $\epsilon_4$ and those obtained with $\epsilon_1$, $\epsilon_2$, and $\epsilon_3$, respectively. The first, second, and third rows display the marginal functions of the joint positions $(\bm x_{15}, \bm x_{16})$ at $t = 0.4 T_f$, $(\bm x_{9}, \bm x_{10})$ at $t = 0.7 T_f$, and $(\bm x_{3}, \bm x_{4})$ at $t = T_f$, respectively. A clear monotonic reduction in the normalized differences is observed as the adaptivity thresholds are tightened, demonstrating convergence of the adaptive TT approximation as higher TT ranks are retained throughout the simulation.

Together with the previous examples, this final demonstration illustrates that the proposed framework can efficiently perform implicit time integration of high-dimensional nonlinear tensor differential equations involving non-diagonal operators and generic non-polynomial nonlinearities, even in regimes where conventional full-order discretizations are computationally intractable.

\section{Conclusions} \label{conclusion}

Tensor train representations have become a powerful tool for solving high-dimensional tensor differential equations arising in scientific computing. While explicit time integration based on cross step-truncation has proven effective for a broad class of problems, extending these methods to implicit schemes has remained an open challenge. The primary obstacle is that the implicit collocation equations associated with the cross-selected fibers are not closed: evaluating the governing operator on a selected fiber requires neighboring fibers that are themselves unknown.

In this work, we introduced a framework that overcomes this obstacle for tensor train representations. The key contribution is a dependency-elimination strategy that expresses the neighboring fibers required by the implicit discretization directly in terms of the \texttt{TT-CROSS-DEIM} interpolation fibers. This converts the implicit collocation equations into a sequence of closed low-dimensional systems that can be solved efficiently without evolving tangent spaces, performing alternating optimization over TT cores, or assembling the full nonlinear tensor. The resulting algorithm preserves the principal advantages of cross step-truncation methods: only $\mathcal{O}(dnr^2)$ tensor entries are evaluated, the computational complexity scales linearly with dimension, and the formulation naturally accommodates generic nonlinearities through pointwise function evaluations.

The proposed framework was demonstrated on both linear and nonlinear tensor differential equations, including three-dimensional and high-dimensional examples with up to thirty spatial dimensions. The numerical results show rapid convergence of the dependency-elimination iteration, preservation of the temporal accuracy of implicit multistep schemes, and effective integration of problems whose corresponding full-order discretizations contain up to $10^{55}$ degrees of freedom. In particular, the results indicate that the convergence of the dependency-elimination iteration is only weakly affected by the size of the underlying tensor discretization or by the presence of non-polynomial nonlinearities.

The proposed methodology establishes implicit cross step-truncation as a practical alternative to existing tensor-train time integration techniques. Beyond the examples considered here, the framework is applicable to a broad class of high-dimensional tensor differential equations arising in uncertainty quantification, kinetic equations, quantum dynamics, stochastic PDEs, and other scientific computing applications involving nonlinear operators.

An interesting observation is that the convergence of the proposed algorithm is governed by the fixed-point iteration used to update the interpolation operators rather than by the conditioning of tangent-space projections or alternating least-squares subproblems. In the numerical examples, the number of fixed-point iterations remains nearly independent of the problem size, suggesting that the contraction properties of the interpolation update may depend only weakly on the ambient dimension. Establishing a rigorous convergence theory for this iteration is an interesting direction for future work.

\appendix
\section{Derivation of the \texttt{TT-CROSS-DEIM} fiber interpolation relations}
\label{appendix:deimrelations}

The purpose of this appendix is to derive
Eqs.~\eqref{eq:deim_dep1}--\eqref{eq:deim_dep3},
which express arbitrary fibers of a \texttt{TT-CROSS-DEIM} approximation as linear
combinations of the DEIM-selected fibers. Throughout this appendix, we use the
interpolatory property of \texttt{TT-CROSS-DEIM},
\[
\hat X(\mathcal I_p,:,\mathcal J_p)
=
X(\mathcal I_p,:,\mathcal J_p),
\]
so that the selected fibers of the original tensor and its \texttt{TT-CROSS-DEIM}
approximation are identical and may be used interchangeably.

For every split \(p=1,\ldots,d-1\), define the left and right contracted
TT vectors
\[
\mathbf g_{1\rightarrow p}
=
G_1(1,i_1,:)\cdots G_p(:,i_p,:)
\in
\mathbb R^{1\times r_p},
\]
and
\[
\mathbf g_{p+1\rightarrow d}
=
G_{p+1}(:,i_{p+1},:)\cdots G_d(:,i_d,1)
\in
\mathbb R^{r_p\times1}.
\]

For the DEIM interpolation index sets
\(\mathcal I_p\) and \(\mathcal J_p\), introduce the interpolation
contraction matrices
\[
\mathbf L_p
=
\begin{bmatrix}
\mathbf g_{1\rightarrow p}^{(1)}
\\
\vdots
\\
\mathbf g_{1\rightarrow p}^{(r_p)}
\end{bmatrix},
\qquad
\mathbf R_p
=
\begin{bmatrix}
\mathbf g_{p+1\rightarrow d}^{(1)}
&
\cdots
&
\mathbf g_{p+1\rightarrow d}^{(r_p)}
\end{bmatrix},
\]
where the superscripts denote evaluation at the interpolation indices.

The following derivation assumes that the selected interpolation fibers have
rank \(r_p\), so that the associated contraction matrices \(\mathbf L_p\) and
\(\mathbf R_p\) are nonsingular. This is the standard interpolation condition
underlying CUR/DEIM approximations. If the prescribed rank exceeds the numerical
rank of the tensor, these matrices may become singular or ill conditioned, in
which case the rank should be reduced or the inverses below interpreted in a
regularized sense. The inverses are introduced only to establish the existence
of interpolation coefficients; the algorithm developed in the main text computes
the coefficient matrices directly from \texttt{TT-CROSS-DEIM} contractions and does not
require forming these inverses explicitly.

We first consider fibers along the first mode. The TT representation gives
\[
\hat X(:,i_2,\ldots,i_d)
=
\mathbf G_1
\mathbf g_{2\rightarrow d},
\]
where \(\mathbf G_1\in\mathbb R^{n_1\times r_1}\) denotes the matricized first
TT core. Evaluating this expression at the selected right interpolation fibers
yields
\[
X(:,\mathcal J_1)
=
\mathbf G_1
\mathbf R_1.
\]
Since \(\mathbf R_1\) is invertible,
\[
\mathbf G_1
=
X(:,\mathcal J_1)\mathbf R_1^{-1},
\]
and therefore
\[
\hat X(:,i_2,\ldots,i_d)
=
X(:,\mathcal J_1)
\mathbf R_1^{-1}
\mathbf g_{2\rightarrow d}.
\]
Thus,
\[
\mathbf z_{r_1}
=
\mathbf R_1^{-1}
\mathbf g_{2\rightarrow d},
\]
which establishes Eq.~\eqref{eq:deim_dep1}.

For an intermediate TT core,
\[
\hat X(i_1,\ldots,i_{p-1},:,i_{p+1},\ldots,i_d)
=
\mathbf g_{1\rightarrow p-1}
\,
G_p
\,
\mathbf g_{p+1\rightarrow d}.
\]
Here and throughout this appendix, products involving a TT core and the interpolation matrices are understood as contractions along the first and third tensor modes. Thus, $\mathbf L_{p-1}
\,G_p\,
\mathbf R_p
$ denotes contraction of \(G_p\) with \(\mathbf L_{p-1}\) along its left TT index and with \(\mathbf R_p\) along its right TT index. Evaluating this expression at the selected interpolation fibers gives
\[
X(\mathcal I_{p-1},:,\mathcal J_p)
=
\mathbf L_{p-1}
\,
G_p
\,
\mathbf R_p.
\]
Hence,
\[
G_p
=
\mathbf L_{p-1}^{-1}
X(\mathcal I_{p-1},:,\mathcal J_p)
\mathbf R_p^{-1},
\]
and substitution into the previous expression yields
\[
\hat X(i_1,\ldots,i_{p-1},:,i_{p+1},\ldots,i_d)
=
\mathbf z_{l_p}
X(\mathcal I_{p-1},:,\mathcal J_p)
\mathbf z_{r_p},
\]
where
\[
\mathbf z_{l_p}
=
\mathbf g_{1\rightarrow p-1}
\mathbf L_{p-1}^{-1},
\qquad
\mathbf z_{r_p}
=
\mathbf R_p^{-1}
\mathbf g_{p+1\rightarrow d}.
\]
This establishes Eq.~\eqref{eq:deim_dep2}.

Finally, for fibers along the last mode,
\[
\hat X(i_1,\ldots,i_{d-1},:)
=
\mathbf g_{1\rightarrow d-1}
G_d.
\]
Evaluating this expression at the selected left interpolation fibers gives
\[
X(\mathcal I_{d-1},:)
=
\mathbf L_{d-1}
G_d,
\]
so that
\[
G_d
=
\mathbf L_{d-1}^{-1}
X(\mathcal I_{d-1},:).
\]
Consequently,
\[
\hat X(i_1,\ldots,i_{d-1},:)
=
\mathbf g_{1\rightarrow d-1}
\mathbf L_{d-1}^{-1}
X(\mathcal I_{d-1},:),
\]
which is Eq.~\eqref{eq:deim_dep3}.

\section{Recursive Identity Property of the \texttt{TT-CROSS-DEIM} Cores}
\label{appendix:identity}

We prove Eq.~\eqref{eq:Gidentity_left} by exploiting the recursive structure of the \texttt{TT-CROSS-DEIM} construction. For the first TT core, the matricized core is constructed as
\[
\mathbf G_1
=
\mathbf U_1
\bigl(
\mathbf U_1(\mathcal I_1,:)
\bigr)^{-1},
\]
where $\mathcal I_1$ denotes the DEIM interpolation indices. Consequently,
\[
\mathbf G_1(\mathcal I_1,:)
=
\mathbf I_{r_1},
\]
which is precisely Eq.~\eqref{eq:Gidentity_left} for $p=1$.

To illustrate the recursive nature of the construction, we next consider the second TT core. After the first TT core has been computed, the remaining tensor is
\[
R_1=X(\mathcal I_1,:)\in
\mathbb R^{r_1\times n_2\times\cdots\times n_d},
\]
where the first mode corresponds to the interpolation index
\(\alpha_1=1,\dots,r_1\). To compute the second TT core, \texttt{CUR-DEIM} is applied to the mode-1 unfolding of \(R_1\), whose rows are indexed by the pair
\[
(\alpha_1,i_2),
\]
obtained by merging the interpolation index \(\alpha_1\) with the physical index \(i_2\). Accordingly,
\[
\mathbf G_2
=
\mathbf U_2
\bigl(
\mathbf U_2(\mathcal I_2,:)
\bigr)^{-1},
\]
where each interpolation index in \(\mathcal I_2\) corresponds to a selected pair
\[
(\alpha_1^{(\alpha_2)},i_2^{(\alpha_2)}),
\qquad
\alpha_2=1,\ldots,r_2.
\]
By construction,
\[
\mathbf G_2(\mathcal I_2,:)
=
\mathbf I_{r_2}.
\]

On the other hand, from the result for the first TT core,
\[
G_1(1,i_1^{(\alpha_1)},:)
=
\mathbf e_{\alpha_1}^{\top},
\]
where \(\mathbf e_{\alpha_1} \in \mathbb R^{r_1 \times 1}\) denotes the \(\alpha_1\)-th canonical basis vector of \(\mathbb R^{r_1}\). Therefore,
\[
G_1(1,i_1^{(\alpha_2)},:)
G_2(:,i_2^{(\alpha_2)},:)
=
G_2(\alpha_1^{(\alpha_2)},i_2^{(\alpha_2)},:),
\]
that is, the contraction of the first two TT cores simply extracts the row of \(\mathbf G_2\) corresponding to the interpolation pair
\((\alpha_1^{(\alpha_2)},i_2^{(\alpha_2)})\). Stacking these rows for all interpolation indices gives
\[
\Big[
G_1(1,i_1^{(\alpha_2)},:)
G_2(:,i_2^{(\alpha_2)},:)
\Big]_{\alpha_2=1}^{r_2}
=
\mathbf G_2(\mathcal I_2,:)
=
\mathbf I_{r_2},
\]
which establishes Eq.~\eqref{eq:Gidentity_left} for \(p=2\).

We now show that the same argument recursively extends to an arbitrary TT core. Assume that
\[
\Big[
G_1(1,i_1^{(\alpha_{p-1})},:)
\cdots
G_{p-1}(:,i_{p-1}^{(\alpha_{p-1})},:)
\Big]_{\alpha_{p-1}=1}^{r_{p-1}}
=
\mathbf I_{r_{p-1}},
\]
and consider the construction of the $p$-th TT core. At this stage, \texttt{TT-CROSS-DEIM} operates on the remainder tensor
\[
R_{p-1}
\in
\mathbb R^{r_{p-1}\times n_p\times\cdots\times n_d}.
\]
Applying \texttt{CUR-DEIM} to its mode-1 unfolding produces
\[
\mathbf G_p
=
\mathbf U_p
\bigl(
\mathbf U_p(\mathcal I_p,:)
\bigr)^{-1},
\]
where the rows of the unfolding are indexed by the merged pair
\[
(\alpha_{p-1},i_p).
\]
Consequently,
\[
\mathbf G_p(\mathcal I_p,:)
=
\mathbf I_{r_p}.
\]

For every interpolation index
\(
(\alpha_{p-1}^{(\alpha_p)},i_p^{(\alpha_p)})
\in\mathcal I_p,
\)
the induction hypothesis implies that
\[
G_1(1,i_1^{(\alpha_p)},:)
\cdots
G_{p-1}(:,i_{p-1}^{(\alpha_p)},:)
=
\mathbf e_{\alpha_{p-1}}^{\top},
\]
where $\mathbf e_{\alpha_{p-1}} \in \mathbb R^{r_{p-1} \times 1}$ is the corresponding canonical basis vector. Therefore,
\[
G_1(1,i_1^{(\alpha_p)},:)
\cdots
G_{p-1}(:,i_{p-1}^{(\alpha_p)},:)
G_p(:,i_p^{(\alpha_p)},:)
=
G_p(\alpha_{p-1}^{(\alpha_p)},i_p^{(\alpha_p)},:),
\]
which is precisely the row of $\mathbf G_p$ selected by the interpolation index. Stacking these rows for all $\alpha_p=1,\ldots,r_p$ gives
\[
\Big[
G_1(1,i_1^{(\alpha_p)},:)
\cdots
G_p(:,i_p^{(\alpha_p)},:)
\Big]_{\alpha_p=1}^{r_p}
=
\mathbf G_p(\mathcal I_p,:)
=
\mathbf I_{r_p},
\]
thereby establishing Eq.~\eqref{eq:Gidentity_left}. Since the argument is independent of $p$, the identity follows recursively for all TT cores.

The proof above establishes the forward interpolation identity
\eqref{eq:Gidentity_left}. The identity
\eqref{eq:Gidentity_right} follows analogously by applying the same
recursive argument to the reverse \texttt{TT-CROSS-DEIM} sweep. In the reverse
construction, the TT cores are computed sequentially from the last mode
toward the first, and the remainder tensors are formed by fixing the
right interpolation indices rather than the left ones. Consequently, the
contracted right TT cores evaluated at the reverse DEIM interpolation
indices satisfy
\[
\mathbf G_{p+1\rightarrow d}
=
\mathbf I_{r_p},
\]
which is the right-to-left counterpart of
Eq.~\eqref{eq:Gidentity_left}. Since the forward and reverse sweeps are
algorithmically identical after reversing the ordering of the tensor
modes, no additional arguments are required.

\section*{Acknowledgments}
This work is sponsored by a grant from Air Force Office of Scientific Research, United States awards no. FA9550-25-1-0039. 

\bibliographystyle{ieeetr} 
\bibliography{References,Hessam}

\end{document}